\documentclass[10pt]{article}

\usepackage{a4wide}
\usepackage{amssymb}
\usepackage{amsfonts}
\usepackage{amsmath}
\input xy
\xyoption{arrow} \xyoption{matrix}

\date{}

\newtheorem{proposition}{Proposition}[section]
\newtheorem{theorem}[proposition]{Theorem}
\newtheorem{lemma}[proposition]{Lemma}

\newtheorem{corollary}[proposition]{Corollary}

\def\Hom{{\rm Hom}}
\def\der{\partial }

\def\nFM0{{\nu }_{F,M_0}}
\def\nFN0{{\nu }_{F,N_0}}
\def\nGN0{{\nu }_{G,N_0}}

\def\N0{ {\bf N}_0 }

\def\t{\otimes}
\def\g{\gamma}
\def\v{\varphi}
\def\ra{\rightarrow}

\def\lra{\leftrightarrow}
\def\Xpm{X^{\pm }}

\def\s{\sigma}
\def\Z{\mathbb{Z}}

\def\l1{{\lambda}_1}

\def\a{\alpha}
\def\a0{ {\alpha }_0}
\def\a1{ {\alpha }_1}

\def\l{\lambda}
\def\o{\omega}

\def\nFGM0{{\nu }_{F,G,M_0}}

\def\nFN0{{\nu}_{F,N_0}}

\def\sm{{\sigma}^m}

\def\sm1{{\sigma}^{-1}}

\def\smtp1{{\sigma}^{-t+1}}

\def\o{\omega }
\def\S1{S^{-1}}

\def\Xpm1{X^{\pm 1}_1}

\def\sPM1{{\sigma }^{\pm 1}}
\def\sMP1{{\sigma }^{\mp 1 }}

\def\d{\delta}

\def\di{{\rm d.ind}}

\def\L{\Lambda}
\def\O{\Omega}

\def\CA{{\cal A}}

\def\CD{{\cal D}}

\def\Ytm1{Y^{t-1}}
\def\Yim1{Y^{i-1}}

\def\CL{{\cal L}}

\def\CG{{\cal G}}
\def\CH{{\cal H}}

\def\Aut{{\rm Aut}}

\def\dim{{\rm dim }}

\def\ker{ {\rm ker } }

\def\CJ{ {\cal J}}

\def\SL2Z{ {\rm SL}_2({\bf Z}) }

\def\CL{{\cal L}}

\def\Gp1{ G^{1 , 1 } }
\def\P11{ P^{-1 , 1 } }
\def\Pp1{ P^{1 , 1 } }

\def\Supp{{\rm Supp}}

\def\nCLsr{{}^\nu\kern-2pt {\cal L}^{\sigma , \rho  }}
\def\nP{{}^\nu \kern-2pt P}
\def\nL{{}^\nu\kern-2pt L}
\def\nLL{{}^\nu\kern-2pt \Lambda}
\def\nPsr{{}^\nu\kern-2pt P^{\sigma , \rho  }}
\def\nLsr{{}^\nu\kern-2pt L^{\sigma , \rho  }}
\def\nuCL{{}^\nu\kern-2pt  {\cal L}}
\def\nCLsr{{}^\nu\kern-2pt {\cal L}^{\sigma , \rho  }}
\def\nCL1m{{}^\nu\kern-2pt {\cal L}^{-1 , 1  }}
\def\x1nu{x^\frac{1}{\nu}}
\def\xm1nu{x^{-\frac{1}{\nu}}}

\def\ra{\rightarrow }

\def\CB{{\cal B}}

\def\CI{{\cal I}}

\def\coker{{\rm coker}}

\def\CC{ {\cal C}}

\def\CH{ {\cal H}}
\def\CP{ {\cal P}}

\def\nAM0{{\nu }_{{\cal A},M_0}}
\def\nAN0{{\nu }_{{\cal A},N_0}}

\def\End{ {\rm End }}

\def\CJ{ {\cal J }}

\def\CP{ {\cal P }}
\def\det{ {\rm det }}

\def\ga{\mathfrak{a}}
\def\gb{\mathfrak{b}}

\def\gd{\mathfrak{d}}

\def\gp{\mathfrak{p}}
\def\gq{\mathfrak{q}}

\def\GL{{\rm GL}}
\def\SL{{\rm SL}}

\def\Spec{{\rm Spec}}

\def\Hom{{\rm Hom}}

\def\di!{\frac{\der^i}{i!}}
\def\dik!{\frac{\der^k_i}{k!}}

\def\id{{\rm id}}

\def\Fun{{\rm Fun}}

\def\N{\mathbb{N}}

\def\0{\overline{0}}
\def\1{\overline{1}}

\def\Ln1{\L_{n,\overline{1}}}

\def\oa{\overline{a}}

\def\a1{a_{\overline{1}}}

\def\St{{\rm St}}
\def\S{\Sigma}

\def\vn1{\overrightarrow{n-1}}

\def\Sh{{\rm Sh}}

\def\im{{\rm im}}

\def\mA{\mathbb{A}}

\def\Sub{{\rm Sub}}
\def\SSub{{\rm SSub}}
\def\Inc{{\rm Inc}}
\def\Min{{\rm Min}}

\def\Inn{{\rm Inn}}

\def\mS{\mathbb{S}}

\def\mI{\mathbb{I}}

\def\lann{{\rm l.ann}}
\def\rann{{\rm r.ann}}
\def\Cen{{\rm Cen}}

\def\hht{{\rm ht}}

\def\mT{\mathbb{T}}
\def\ind{{\rm ind}}

\def\mG{\mathbb{G}}
\def\mE{\mathbb{E}}

\def\mU{\mathbb{U}}

\def\Out{{\rm Out}}

\def\rG{\mathrm{G}}
\def\brG{\overline{\rG}}
\def\Supp{{\rm Supp}}

\begin{document}

\author{V. V. \  Bavula 
}

\title{The group of automorphisms of the   algebra of polynomial  integro-differential operators}

\maketitle

\begin{abstract}
The group $\rG_n$ of automorphisms of  the   algebra
$\mI_n:=K\langle x_1, \ldots , x_n, \frac{\der}{\der x_1}, \ldots
,\frac{\der}{\der x_n}, \int_1, \ldots , \int_n\rangle $ of
polynomial  integro-differential operators  is  found:
$$ \rG_n=S_n\ltimes \mT^n\ltimes  \Inn (\mI_n) \supseteq
 S_n\ltimes \mT^n \ltimes \underbrace{\GL_\infty (K)\ltimes\cdots \ltimes
\GL_\infty (K)}_{2^n-1 \;\; {\rm times}}, $$
$$ \rG_1\simeq  \mT^1 \ltimes \GL_\infty (K),$$
 where  $S_n$
is the symmetric group, $\mT^n$ is the $n$-dimensional torus,
  $\Inn (\mI_n)$ is
the group of inner automorphisms of $\mI_n$ (which is huge). It is
proved that each automorphism $\s \in \rG_n$ is uniquely
determined by the elements $\s (x_i)$'s or $\s (\frac{\der}{\der
x_i})$'s or $\s (\int_i)$'s.
 The stabilizers in $\rG_n$ of all the ideals of
$\mI_n$ are found, they are subgroups of {\em finite}  index in
$\rG_n$. It is shown that the group $\rG_n$ has trivial centre,
$\mI_n^{\rG_n}=K$ and $\mI_n^{\Inn (\mI_n)}=K$, the (unique)
maximal ideal of $\mI_n$ is the {\em only} nonzero prime
$\rG_n$-invariant ideal of $\mI_n$, and there are precisely $n+2$
$\rG_n$-invariant ideals of $\mI_n$.
 For each automorphism $\s \in \rG_n$, an {\em explicit inversion
 formula}
is given  via the elements $\s (\frac{\der}{\der x_i})$ and  $\s
(\int_i)$.



$${\bf Contents}$$
\begin{enumerate}
\item Introduction. \item The algebras  $\mI_n$ and $\mA_n$. \item
A description of the group $\rG_n$ and two criteria.  \item The
group $\Aut_{K-{\rm alg}}(\mI_1)$. \item The group of
automorphisms of the algebra $\mI_n$. \item Stabilizers of the
ideals of $\mI_n$ in $\rG_n$.
\end{enumerate}
\end{abstract}


\section{Introduction}\label{INTRO}
Throughout, ring means an associative ring with $1$; module means
a left module;
 $\N :=\{0, 1, \ldots \}$ is the set of natural numbers; $K$ is a
field of characteristic zero and  $K^*$ is its group of units;
$P_n:= K[x_1, \ldots , x_n]$ is a polynomial algebra over $K$;
$\der_1:=\frac{\der}{\der x_1}, \ldots , \der_n:=\frac{\der}{\der
x_n}$ are the partial derivatives ($K$-linear derivations) of
$P_n$; $\End_K(P_n)$ is the algebra of all $K$-linear maps from
$P_n$ to $P_n$ and $\Aut_K(P_n)$ is its group of units (i.e. the
group of all the invertible linear maps from $P_n$ to $P_n$); the
subalgebra  $A_n:= K \langle x_1, \ldots , x_n , \der_1, \ldots ,
\der_n\rangle$ of $\End_K(P_n)$ is called the $n$'th {\em Weyl}
algebra.

$\noindent $

{\it Definition}, \cite{Bav-Jacalg}. The {\em Jacobian algebra}
$\mA_n$ is the subalgebra of $\End_K(P_n)$ generated by the Weyl
algebra $A_n$ and the elements $H_1^{-1}, \ldots , H_n^{-1}\in
\End_K(P_n)$ where $$H_1:= \der_1x_1, \ldots , H_n:= \der_nx_n.$$

Clearly, $\mA_n =\bigotimes_{i=1}^n \mA_1(i) \simeq \mA_1^{\t n }$
where $\mA_1(i) := K\langle x_i, \der_i , H_i^{-1}\rangle \simeq
\mA_1$. The algebra $\mA_n$ contains all the  integrations
$\int_i: P_n\ra P_n$, $ p\mapsto \int p \, dx_i$, since  $$\int_i=
x_iH_i^{-1}: x^\alpha \mapsto (\alpha_i+1)^{-1}x_ix^\alpha.$$ In
particular, the algebra $\mA_n$ contains the {\em algebra}
$\mI_n:=K\langle x_1, \ldots , x_n$,
 $\der_1, \ldots ,\der_n,  \int_1,
\ldots , \int_n\rangle $ {\em of  polynomial integro-differential
operators}. Note that $\mI_n=\bigotimes_{i=1}^n\mI_1(i)\simeq
\mI_1^{\t n}$ where $\mI_1(i):= K\langle x_i, \der_i,
\int_i\rangle$. Let $\rG_n:=\Aut_{K-{\rm alg}}(\mI_n)$ and $\mG_n
:=\Aut_{K-{\rm alg}}(\mA_n)$.

The Jacobian algebra $\mA_n$ is a (two-sided) localization $\mA_n
= S^{-1}\mI_n$ of the algebra $\mI_n$  at a countably generated
commutative monoid $S\simeq \N^{(\N )}$, each element of $S$ is a
regular element of the algebra $\mI_n$, \cite{algintdif}. In
general, there is no connection between the groups of
automorphisms of an algebra and its localization. As a rule,  the
latter is smaller than the former, and an automorphism of the
algebra cannot be extended  to an automorphism of its localization
(eg, the group $\Aut_{K-{\rm alg}}(P_n)$ is huge but $\Aut_{K-{\rm
alg}}(K[x_1^{\pm 1}, \ldots , x_n^{\pm 1}]$ is a tiny group).
Completely the opposite is true for  the pair $\mI_n$, $\mA_n =
S^{-1}\mI_n$: each automorphism of the algebra $\mI_n$ can be
extended to an automorphism of the algebra $\mA_n$ (this is not
straightforward since $\rG_n S\not\subseteq S$). Moreover, the
group $\rG_n$ can be seen as a subgroup of $\mG_n$ (Theorem
\ref{B25Oct9}), and the group $\mG_n$ is bigger than $\rG_n$. This
fact, i.e. $\rG_n\subseteq \mG_n$, is one of the key moments in
finding the group $\rG_n$ as the group $\mG_n$ was already found
in \cite{jacaut}.

The algebras $P_{2n}=P_2^{\t n}$, $A_n = A_1^{\t n}$, $\mS_n =
\mS_1^{\t n}$, $\mI_n = \mI_1^{\t n}$ and $\mA_n = \mA_1^{\t n}$
have similar defining relations:
\begin{eqnarray*}
 P_2=K\langle x,y \rangle &:& yx-xy=0;\\
A_1=K\langle x, \der\rangle &:& \der x-x\der =1;\\
\mS_1=K\langle x, y\rangle &:& yx=1;\\
\mI_1=K\langle\der, H, \int  \rangle &:& \der \int =1, [ H, \int ]
= \int, [ H, \der ] = -\der,
H(1-\int\der ) = (1-\int \der ) H = 1-\int \der ; \\
\mA_1=K\langle x,H^{\pm 1}, y \rangle &:& yx=1, [H, x]=x, [H,
y]=-y, H(1-xy) = (1-xy)H= 1-xy;
\end{eqnarray*}
where  $[a,b] := ab-ba$ is the commutator of elements $a$ and $b$.
It is reasonable to believe that that they should have similar
groups of automorphisms. This is exactly  the case when $n=1$: the
groups of automorphisms of the algebras $P_2$, $A_1$, $\mS_1$,
$\mI_1$ and $\mA_1$ have almost identical structure (when properly
interpreted). Namely, each of the groups is a `product' (in the
last three cases  it is even the semi-direct product) of an
obvious subgroup of affine automorphisms and a non-obvious
subgroup generated by `transvections.'

The group $\Aut_{K-{\rm alg}}(P_2)$ was found by Jung \cite{jung}
in 1942 and van der Kulk \cite{kulk} in 1953. In 1968, Dixmier
\cite{Dix} found the group of automorphisms of the first Weyl
algebra $A_1$ (in prime characteristic  the group of automorphism
of the first Weyl algebra $A_1$ was found by Makar-Limanov
\cite{Mak-LimBSMF84} in 1984, see also \cite{A1rescen}  for a
different approach and for further developments). In 2000,
Gerritzen \cite{Gerritzen-2000} found generators for the group
$\Aut_{K-{\rm alg}}(\mS_1)$.
 For
the higher Weyl algebras $A_n$, $n\geq 2$, and the polynomial
algebras $P_n$, $n\geq 3$,  to find their groups of automorphisms
and generators are old open problems, and solutions to the
Jacobian Conjecture and the Problem/Conjecture of Dixmier would be
an important (easier) part in finding the groups (positive
solutions to these two problems would define the groups as
infinite dimensional varieties, i.e. they would  give defining
equations of the varieties but not generators. To find generators
one would have to find the solutions of the equations. A finite
dimensional analogue of this situation is the group $\SL_n$: the
defining equation $\det =1$ tells nothing about generators of the
group).

The Jacobian algebras $\mA_n$  arose in my study of the group of
polynomial automorphisms and the Jacobian Conjecture, which is a
conjecture that makes sense {\em only} for polynomial algebras in
the class of all commutative algebras \cite{Bav-inform}. In order
to solve the Jacobian Conjecture,   it is reasonable to believe
that one should create a technique which makes sense {\em only}
for polynomials;  the Jacobian algebras are a step in this
direction (they exist for polynomials but make no sense even for
Laurent polynomials).

The Jacobian algebras $\mA_n$ were invented to deal with
polynomial automorphisms. A study of these algebras led to study
of `simpler' algebras $\mS_n$ \cite{shrekalg}, the so-called {\em
algebras of one-sided inverses of  polynomial algebras}. This
ended up in finding their groups of automorphisms $\Aut_{K-{\rm
alg}}(\mS_n)$, $n\geq 1$, and their explicit generators in the
series of three papers \cite{shrekaut}, \cite{K1aut} and
\cite{Snaut}. Recently, the groups $\Aut_{K-{\rm alg}}(\mA_n)$,
$n\geq 1$, are found in \cite{jacaut}. Finally, in the present
paper the groups $\rG_n:=\Aut_{K-{\rm alg}}(\mI_n)$, $n\geq 1$,
are found.
\begin{itemize}
\item (Theorem \ref{25Oct9}.(1)) $\rG_n = S_n\ltimes \mT^n \ltimes
\Inn (\mI_n)$ {\em  where $S_n$ is the symmetric group, $\mT^n$ is
the $n$-dimensional torus and $\Inn (\mI_n)$ is the group of inner
automorphisms of the algebra} $\mI_n$. \item (Theorem
\ref{A25Oct9}.(2)) {\em The map $(1+\ga_n)^*\ra \Inn (\mI_n)$,
$u\mapsto \o_u$, is a group isomorphism where $\o_u (a) :=
uau^{-1}$, $(1+\ga_n)^* := \mI_n^* \cap (1+\ga_n)$, $\mI_n^*$ is
the group of units of the algebra $\mI_n$, and $\ga_n$ is the only
maximal ideal of the algebra $\mI_n$.}
\end{itemize}

$\noindent $

The paper proceeds as follows. In Section \ref{TAIAN}, some known
results about the algebras $\mI_n$ and $\mA_n$ are collected that
are used freely in the paper.

One of the key ideas in finding the group $\rG_n$ is the fact that
the polynomial algebra $P_n$ is the only (up to isomorphism)
faithful simple $\mI_n$-module (Proposition 3.8,
\cite{algintdif}). This enables us to visualize the group $\rG_n$
as a subgroup of $\Aut_K(P_n)$ (Corollary \ref{d21Mar9}):
$$ \rG_n = \{ \s_\v \, | \, \v \in \Aut_K(P_n), \; \v \mI_n
\v^{-1} = \mI_n\} \;\; {\rm where}\;\; \s_\v (a) :=\v a \v^{-1},
\; a\in \mI_n.$$ In Section \ref{ADGN}, two `rigidity theorems'
are proved for the group $\rG_n$: Corollary \ref{yz16Apr9} and
\begin{itemize}
\item (Theorem \ref{21Mar9}) {\rm (Rigidity  of the group
$\rG_n$)} {\em Let $\s , \tau \in \rG_n$. Then the following
statements are equivalent.}
\begin{enumerate}
 \item $\s = \tau$. \item $\s (\int_1) = \tau (\int_1), \ldots , \s (\int_n) = \tau
 (\int_n)$.
 \item $\s (\der_1) = \tau (\der_1), \ldots , \s (\der_n) = \tau (\der_n)$.
 \item $\s (x_1) = \tau (x_1), \ldots , \s (x_n) = \tau (x_n)$.
\end{enumerate}
\end{itemize}
In Section \ref{AAII1}, the group $\rG_1$ and its explicit
generators are found (Theorem \ref{31Oct9}). The key ingredients
of the proof of Theorem \ref{31Oct9} are Fredholm operators, their
indices and the Rigidity of the group $\rG_1$ (Theorem
\ref{21Mar9}). It is proved that each algebra endomorphism of the
algebra $\mI_1$ is a monomorphism (Theorem \ref{B31Oct9}), and no
proper prime factor algebra of the algebra $\mI_n$ can be embedded
into the algebra $\mI_n$ (Theorem \ref{B14Nov9}). These two
results have bearing of the Jacobian Conjecture and the
Problem/Conjecture of Dixmier ({\em each algebra endomorphism of
the Weyl algebra is an isomorphism}).

Section \ref{GAAIN} contains the main results of the paper, a
proof that $\rG_n = S_n\ltimes \mT^n \ltimes \Inn (\mI_n)$
(Theorem \ref{25Oct9}) is given.

\begin{itemize}
\item (Theorem \ref{B25Oct9}) \begin{enumerate} \item $\rG_n=\{ \s
\in \mG_n \, | \, \s (\mI_n) = \mI_n\}$ {\em and $\rG_n$ is a
subgroup of} $\mG_n$. \item {\em Each automorphism of the algebra
$\mI_n$ has a unique extension to an automorphism of the algebra}
$\mI_n$.
\end{enumerate}
\item (Theorem \ref{14Nov9}) {\em The centre of the group $\rG_n$
is} $\{ e \}$. \item (Theorem \ref{A14Nov9}) $\mI_n^{\rG_n}=K$
{\em and $ \mI_n^{\Inn (\mI_n)}=K$, the  algebras of invariants.}
\end{itemize}

Each automorphism $\s \in \rG_n =S_n\ltimes \mT^n\ltimes \Inn
(\mI_n)$ is a unique product $\s = st_\l \o_\v$ which is called
the {\em canonical form} of $\s$ where $s\in S_n$, $t_\l \in
\mT^n$, $\o_\v \in \Inn (\mI_n)$ and $\v \in (1+\ga_n)^*$ ($\v $
is unique).

\begin{itemize}
\item (Corollary \ref{a2Nov9}) {\em Let $\s \in \rG_n$ and $ \s
=st_\l\o_\v$ be its canonical form. Then the automorphisms $s$,
$t_\l$,  and $\o_\v$ can be effectively (in finitely many steps)
found from the action of the automorphism $\s$ on the elements $\{
H_i, \der_i, \int_i\, | \, i=1, \ldots , n\}$:
$$ \s (H_i) \equiv H_{s(i)}\mod \ga_n, \;\; \s (\der_i) \equiv
\l_i^{-1} \der_{s_(i)}\mod \ga_n,\;\; \s (\int_i) \equiv
\l_i\int_{s(i)}\mod \ga_n,
$$ and the elements $\v$  and $\v^{-1}$ are  given by the formulae
 (\ref{vf2Nov9}) and (\ref{1vf2Nov9}) respectively
for the automorphism} $(st_\l)^{-1}\s \in \Inn (\mI_n)$.
\end{itemize}
The explicit formulae (\ref{1vf2Nov9}) and (\ref{vf2Nov9}) are too
complicated to reproduce them in the Introduction.
\begin{itemize}
\item (Corollary \ref{b2Nov9}) (A criterion of being inner
automorphism)  {\em Let $\s \in \rG_n$. The following statements
are equivalent.}
\begin{enumerate}
\item $\s \in \Inn (\mI_n)$. \item $\s (\der_i) \equiv \der_i \mod
\ga_n$ {\em for} $i=1, \ldots , n$. \item $\s (\int_i) \equiv
\int_i \mod \ga_n$ {\em for} $i=1, \ldots , n$.
\end{enumerate}
\end{itemize}
{\bf An inversion formula for $\s\in \rG_n$}.
 The next theorem
gives an inversion formula for $\s\in \rG_n$ via the elements $\{
\s (\der_i) , \s(\int_i)\, | \, i=1, \ldots , n\}$.
\begin{itemize}
\item (Theorem \ref{13Nov9}) {\em Let $\s \in \rG_n$ and $\s =
st_\l \o_\v$ be its canonical form where $s\in S_n$, $t_\l\in
\mT^n$ and $\o_\v \in \Inn (\mI_n)$ for a unique element $\v \in
(1+\ga_n)^*$. Then $\s^{-1} = s^{-1} t_{s(\l^{-1})}\o_{st_\l
(\v^{-1})}$  is the canonical form of the automorphism $\s^{-1}$
where the elements  $\v^{-1}$ and $\v$  are given by
 the formulae (\ref{1vf2Nov9}) and (\ref{vf2Nov9}) respectively for the
 automorphism}
 $(st_\l )^{-1}\s \in \Inn (\mI_n)$.
\end{itemize}

In Section \ref{ISTABAN}, the stabilizers in the group $\rG_n$ of
 all the ideals of the algebra $\mI_n$ are computed (Theorem \ref{15Apr9}).
 In particular, the
stabilizers of all the prime ideals of $\mI_n$ are found
(Corollary \ref{b15Apr9}.(2)). 

\begin{itemize}
\item (Corollary \ref{b15Apr9}.(3)) {\em The ideal $\ga_n$ is the
only nonzero, prime, $\rG_n$-invariant ideal of the algebra
$\mI_n$.} \item  (Corollary  \ref{b15Apr9}) {\em Let $\gp$ be a
prime ideal of $\mI_n$. Then its stabilizer $\St_{\mI_n}(\gp )$ is
a maximal subgroup of the group $\rG_n$ iff $n>1$ and $\gp$ is of
height 1, and, in this case, $[\rG_n:\St_{\rG_n}(\gp )]=n$. }\item
(Corollary \ref{d15Apr9}) {\em Let $\ga$ be a proper ideal of
$\mI_n$. Then its stabilizer $\St_{\rG_n}(\ga )$ has finite index
in the group $\rG_n$.} \item (Corollary \ref{c15Apr9}) {\em If
$\ga$ is a generic  ideal of $\mI_n$ then its stabilizer can be
written via the wreath products of the symmetric groups:}
$$  \St_{\rG_n}(\ga )= (S_m\times \prod_{i=1}^t(S_{h_i}\wr
S_{n_i}))\ltimes \mT^n\ltimes \Inn (\mI_n),  $$ {\em where $\wr$
stands for the {\em wreath} product of groups.}
\end{itemize}
Corollary \ref{e15Apr9} classifies all the proper
$\rG_n$-invariant ideals of the algebra $\mI_n$, there are exactly
$n$ of them.


\section{The algebras  $\mI_n$ and $\mA_n$}\label{TAIAN}

In this section, for the reader's convenience we collect some
known results about the algebras $\mI_n$ and $\mA_n$ from the
papers \cite{Bav-Jacalg, jacaut, algintdif} that are used later in
the paper.

The algebra $\mI_n$ is a prime, central, catenary, non-Noetherian
algebra of classical Krull dimension $n$ and of Gelfand-Kirillov
dimension $2n$, \cite{algintdif}. Since $x_i= \int_iH_i$, the
algebra $\mI_n$ is generated by the elements $\{ \der_i , H_i,
\int_i\, | \, i=1, \ldots , n\}$, and $\mI_n =\bigotimes_{i=1}^n
\mI_1(i)$ where $\mI_1(i):= K\langle \der_i , H_i,
\int_i\rangle=K\langle \der_i , x_i, \int_i\rangle\simeq \mI_1$.
 When $n=1$ we usually drop the subscript `1' in $\der_1$, $\int_1$, $H_1$, and $x_1$.
 The following elements of the algebra $\mI_1=K\langle
\der , H, \int \rangle$, 
\begin{equation}\label{eijdef}
e_{ij}:=\int^i\der^j-\int^{i+1}\der^{j+1}, \;\; i,j\in \N ,
\end{equation}
satisfy the relations: $e_{ij}e_{kl}=\d_{jk}e_{il}$ where
$\d_{ij}$ is the Kronecker delta. The matrices of the linear maps
$e_{ij}\in \End_K(K[x])$  with respect to the basis $\{ x^{[s]}:=
\frac{x^s}{s!}\}_{s\in \N}$ of the polynomial algebra $K[x]$  are
the elementary matrices, i.e.
$$ e_{ij}*x^{[s]}=\begin{cases}
x^{[i]}& \text{if }j=s,\\
0& \text{if }j\neq s.\\
\end{cases}$$
The direct sum  $F:=\bigoplus_{i,j\in \N}Ke_{ij}$ is the only
proper (hence maximal) ideal of the algebra $\mI_1$. As an algebra
without 1 it is isomorphic to the algebra without 1 of infinite
dimensional matrices $M_\infty (K) :=\varinjlim
M_d(K)=\bigoplus_{i,j\in \N}KE_{ij}$ via $e_{ij}\mapsto E_{ij}$
where $E_{ij}$ are the matrix units. For all $i,j\in \N$,
\begin{equation}\label{Ieij}
\int e_{ij}=e_{i+1,j}, \;\; e_{ij}\int = e_{i,j-1}, \;\;\der
e_{ij}=e_{i-1,j}, \;\; e_{ij}\der= e_{i,j+1},
\end{equation}
where $e_{-1, j}:=0$ and $e_{i,-1}:=0$. The algebra
$\mI_n=\bigotimes_{i=1}^n \mI_1(i)$ contains the ideal $F_n:=
F^{\t n}=\bigotimes_{i=1}^nF(i)=\bigoplus_{\alpha, \beta \in \N^n}
Ke_{\alpha\beta}$ where $e_{\alpha\beta}:= \prod_{i=1}^n
e_{\alpha_i\beta_i}(i)$,
$e_{\alpha_i\beta_i}(i):=\int_i^{\alpha_i}\der_i^{\beta_i}-
\int_i^{\alpha_i+1}\der_i^{\beta_i+1}$ and $F(i)=\bigoplus_{s,t\in
\N}Ke_{st}(i)$.

\begin{proposition}\label{a5Oct9}
{\rm \cite{algintdif}}
\begin{enumerate}
\item The algebra $\mI_n$ is generated by the elements $\{ \der_i,
\int_i, H_i\, | \, i=1, \ldots , n\}$ that satisfy the following
defining relations:
\begin{eqnarray*}
\forall i& : & \der_i\int_i= 1, \;\; [H_i, \int_i]=\int_i, \;\;
[ H_i, \der_i]=-\der_i, \;\; H_i (1-\int_i\der_i) = (1-\int_i\der_i) H_i=1-\int_i\der_i, \\
\forall i\neq j &:   & a_ia_j= a_ja_i \;\; {\rm where}\;\; a_k\in \{ \der_k, \int_k, H_k\}. \\
\end{eqnarray*}
\item The algebra $\mI_n= \bigotimes_{i=1}^n D_1(i)(\s_i,
1)=D_n((\s_1, \ldots , \s_n), (1,\ldots , 1))$ is a generalized
Weyl algebra ($\int_i\lra x_i$, $\der_i\lra y_i$, $H_i\lra H_i$)
where $D_n:=\bigotimes_{i=1}^n D_1(i)$, $D_1(i):=K[H_i]\bigoplus
\bigoplus_{j\geq 0}Ke_{jj}(i)$, $H_ie_{jj}(i) =
e_{jj}(i)H_i=(j+1)e_{jj}(i)$,  and the $K$-algebra endomorphisms
$\s_i$ are given by the rule $\s_i (a):= \int_ia\der_i$
($\s_i(H_i) = H_i-1$, $\s_i(e_{jj}(i))=e_{j+1, j+1}(i)$).
Moreover, the algebra $\mI_n=\bigoplus_{\alpha \in \Z^n}\mI_{n,
\alpha}$ is $\Z^n$-graded where $\mI_{n,\alpha}= D_n v_\alpha =
v_\alpha D_n $ for all $\alpha \in \Z^n$ where $v_\alpha:=
\prod_{i=1}^nv_{\alpha_i}(i)$ and $v_j(i):=\begin{cases}
\int_i^j& \text{if } j>0,\\
1& \text{if } j=0,\\
\der_i^{-j} & \text{if } j<0.\\
\end{cases} $
\end{enumerate}
\end{proposition}
{\it Remark}. Note that $\s_i(1) = \int_i\der_i = 1-e_{00}(i)\neq
1$ for all $i=1, \ldots , n$.

$\noindent $

 {\it Definition}. Let $A$ and $B$
be algebras, let $\CJ (A)$ and $\CJ (B)$ be their lattices of
ideals. We say that a bijection $f: \CJ (A) \ra \CJ (B)$ is an
{\em isomorphism} if $f(\ga *\gb ) = f(\ga )*f(\gb )$ for $*\in \{
+, \cdot, \cap \}$, and in this case we say that the algebras $A$
and $B$ are {\em ideal equivalent}.

$\noindent $

The ideal equivalence is an equivalence relation on the class of
algebras (introduced in \cite{algintdif}). The next theorem shows
that the Jacobian algebra $\mA_n$ and the algebra $\mI_n$ are
ideal equivalent.
\begin{theorem}\label{7Oct9}
{\em \cite{algintdif}} The restriction map $\CJ ( \mA_n) \ra \CJ
(\mI_n)$, $\ga \mapsto \ga^r:= \ga \cap \mI_n$, is an isomorphism
(i.e. $(\ga_1*\ga_2)^r= \ga_1^r*\ga_2^r$ for $*\in \{ +, \cdot,
\cap \}$) and its inverse is the extension map $\gb\mapsto \gb^e:=
\mA_n \gb \mA_n$.
\end{theorem}

 The next corollary shows that the ideal
theory of $\mI_n$ is `very arithmetic.' In some sense, it is the
best and the simplest possible ideal theory one can imagine.  Let
$\CB_n$ be the set of all functions $f:\{ 1, 2, \ldots , n\} \ra
\{ 0,1\}$. For each function $f\in \CB_n$, $I_f:= I_{f(1)}\t
\cdots \t I_{f(n)}$ is the ideal of $\mI_n$ where $I_0:=F$ and
$I_1:= \mI_1$.  Let $\CC_n$ be the set of all subsets of $\CB_n$
all distinct elements of which are incomparable (two distinct
elements $f$ and $g$ of $\CB_n$ are {\em incomparable} if neither
$f(i)\leq g(i)$ nor $f(i)\geq g(i)$ for all $i$). For each $C\in
\CC_n$, let $I_C:= \sum_{f\in C}I_f$, the ideal of $\mI_n$. The
number $\gd_n$ of elements in the set $\CC_n$ is called the {\em
Dedekind number}. It appeared in the paper of Dedekind
\cite{Dedekind-1871}. An asymptotic of the Dedekind numbers was
found by Korshunov \cite{Korshunov-1977}.

Recall that a submodule of a module that intersects non-trivially
each nonzero submodule of the module is called an {\em essential}
submodule.

\begin{corollary}\label{b10Oct9}
{\em \cite{algintdif}}
\begin{enumerate}
\item The set of height one prime ideals of the algebra $\mI_n$ is
$\{ \gp_1:= F\t\mI_{n-1}, \gp_1:= \mI_1\t F\t\mI_{n-2},\ldots ,
\gp_n:= \mI_{n-1}\t F\}$. \item Each ideal of the algebra $\mI_n$
is an idempotent ideal ($\ga^2= \ga$). \item The ideals of the
algebra $\mI_n$ commute ($\ga \gb = \gb \ga$). \item The lattice
$\CJ (\mI_n)$ of ideals of the algebra $\mI_n$ is distributive.
\item Each nonzero ideal of the algebra $\mI_n$ is an essential
left and right submodule of $\mI_n$. \item $\ga \gb = \ga \cap
\gb$ for all ideals $\ga $ and $\gb $ of the algebra $\mI_n$.
\item The ideal $\ga_n := \gp_1+\cdots + \gp_n$ is the largest
(hence, the only maximal) ideal of $\mI_n$ distinct from $\mI_n$,
and $F_n = F^{\t n}=\cap_{i=1}^n \gp_i$ is the smallest nonzero
ideal of $\mI_n$. \item {\rm (A classification of ideals of
$\mI_n$)} The map
 $\CC_n\ra \CJ (\mI_n)$, $C\mapsto I_C:= \sum_{f\in C}I_f$
 is a bijection where $I_\emptyset :=0$. The number of ideals of
 $\mI_n$ is equal to the Dedekind  number $\gd_n$.
 For $n=1$, $F$ is the unique proper ideal of the algebra
$\mI_1$.  \item {\rm (A classification of prime ideals of
$\mI_n$)} Let $\Sub_n$ be the set of all subsets of $\{ 1, \ldots
, n\}$. The map $\Sub_n\ra \Spec (\mI_n)$, $ I\mapsto \gp_I:=
\sum_{i\in I}\gp_i$, $\emptyset \mapsto 0$, is a bijection, i.e.
any nonzero prime ideal of $\mI_n$ is a unique sum of primes of
height 1; $|\Spec (\mI_n)|=2^n$; the height of $\gp_I$ is $| I|$;
and
 \item  $\gp_I\subset \gp_J$ iff $I\subset
 J$.
\end{enumerate}
\end{corollary}

{\bf The involution $*$ on the algebra $\mI_n$}. Using the
defining relations in Proposition \ref{a5Oct9}.(1), we see that
the algebra $\mI_n$ admits the involution: 
\begin{equation}\label{*invIn}
*: \mI_n\ra \mI_n, \;\; \der_i \mapsto  \int_i, \;\; \int_i\mapsto
\der_i, \;\; H_i\mapsto H_i, \;\; i=1, \ldots , n,
\end{equation}
i.e. it is a $K$-algebra {\em anti-isomorphism} $((ab)^*= b^*a^*)$
such that $*\circ *= \id_{\mI_n}$. Therefore, the algebra $\mI_n$
is {\em self-dual}, i.e. is isomorphic to its {\em opposite}
algebra $\mI_n^{op}$. As a result, the left and the right
properties of the algebra $\mI_n$ are the same.  For all elements
$\alpha , \beta \in \N^n$, 
\begin{equation}\label{eab*}
e_{\alpha\beta}^*= e_{\beta \alpha}.
\end{equation}

The involution $*$ can be extended to an involution of the algebra
$\mA_n$ by setting
$$ x_i^*=H_i\der_i, \;\; \der_i^* = \int_i, \;\; (H_i^{\pm 1})^*=
H_i^{\pm 1}, \;\; i=1, \ldots , n.$$  Note that $y_i^* =
(H_i^{-1}\der_i)^* = \int_i H_i^{-1} = x_iH_i^{-2}$, $A_n^*
\not\subseteq A_n$,  but $\CI_n^* = \CI_n$ where $\CI_n:= K\langle
\der_1, \ldots , \der_n \int_1, \ldots , \int_n\rangle $ is the
algebra of integro-differential operators with constant
coefficients.

For a subset $S$ of a ring $R$, the sets $\lann_R(S):= \{ r\in R\,
| \, rS=0 \}$ and $\rann_R(S):= \{ r\in R\, | \, Sr=0\}$ are
called the {\em left} and the {\em right annihilators} of the set
$S$ in $R$.
 Using the fact that the algebra $\mI_n$ is a GWA and its
 $\Z^n$-grading, we see that
\begin{equation}\label{laI}
\lann_{\mI_n}(\int_i) = \bigoplus_{k\in \N}Ke_{k0}(i)\bigotimes
\bigotimes_{i\neq j}\mI_1(j), \;\; \rann_{\mI_n}(\int_i) =0.
\end{equation}
\begin{equation}\label{lad}
\rann_{\mI_n}(\der_i) = \bigoplus_{k\in \N}Ke_{0k}(i)\bigotimes
\bigotimes_{i\neq j}\mI_1(j), \;\; \lann_{\mI_n}(\der_i) =0.
\end{equation}

Let $\ga$ be an ideal of the algebra $\mI_n$. The factor algebra
$\mI_n / \ga$ is a Noetherian algebra iff $\ga = \ga_n$,
\cite{algintdif}. The factor algebra $B_n := \mI_n / \ga_n$ is
isomorphic to the skew Laurent polynomial algebra
$$\bigotimes_{i=1}^n K[H_i][\der_i,
\der_i^{-1}; \tau_i]=\CP_n [ \der_1^{\pm 1}, \ldots, \der_n^{\pm
1}; \tau_1, \ldots , \tau_n],$$ via $\der_i\mapsto \der_i$,
$\int_i\mapsto \der_i^{-1}$, $H_1\mapsto H_i$  (and $x_i\mapsto
\der_i^{-1}H_i$) where $\CP_n := K[H_1, \ldots , H_n]$ and
$\tau_i(H_i) = H_i+1$. We identify these two algebras via this
isomorphism. It is obvious that
$$ B_n =\bigotimes_{i=1}^n K[H_i][z_i,
z_i^{-1}; \s_i]=\CP_n [ z_1^{\pm 1}, \ldots, z_n^{\pm 1}; \s_1,
\ldots , \s_n],$$ where $z_i:=\der_i^{-1}$ and  $\s_i =
\tau_i^{-1}: H_i\mapsto H_i-1$. We use this alternative
presentation of the algebra $B_n$ in order to avoid awkward
expressions like $\frac{\der}{\der \der_i}$ later.
 By Theorem \ref{7Oct9}, $\ga_n^e$  is the only maximal ideal of
 the Jacobian algebra $\mA_n$. The factor algebra $\CA_n := \mA_n
 / \ga_n^e$ is the skew Laurent polynomial algebra
 $$ \CA_n =\CL_n [\der_1^{\pm 1}, \ldots , \der_n^{\pm 1}; \tau_1, \ldots
 , \tau_n] = \CL_n [x_1^{\pm 1}, \ldots , x_n^{\pm 1}; \s_1, \ldots
 , \s_n] = \CL_n [z_1^{\pm 1}, \ldots , z_n^{\pm 1}; \s_1, \ldots
 , \s_n]$$ where $\CL_n:= K[H_1^{\pm 1}, (H_1\pm 1)^{-1}, (H_1\pm 2)^{-1},
\ldots , H_n^{\pm 1}, (H_n\pm 1)^{-1}, (H_n\pm 2)^{-1}, \ldots ]$,
$\tau_i(H_j)= H_j+\d_{ij}$ where $\d_{ij}$ is the Kronecker delta
and  $\s_i=\tau_i^{-1}$. By Theorem \ref{7Oct9},
$\ga_n^{er}=\ga_n$, hence  the algebra $B_n$ is a subalgebra of
$\CA_n$. Moreover, the algebra $\CA_n$ is the localization of the
algebra $B_n$ at the multiplicatively closed set $\{
(H_1+\alpha_1)^{m_1}\cdots (H_n +\alpha_n)^{m_n}\, | \,
(\alpha_i)\in \Z^n , (m_i) \in \N^n\}$. The algebra $B_n$ is also
the left (but not right) localization of the algebra $\mI_n$ at
the multiplicatively closed set $S_{\der_1, \ldots , \der_n}:=\{
\der_1^{\alpha_1}\cdots \der_n^{\alpha_n}\, | \, (\alpha_i)\in
\N^n\}$, $B_n \simeq S_{\der_1, \ldots , \der_n}^{-1}\mI_n$.


\section{A description of the group $\rG_n$ and two criteria}\label{ADGN}

In this section, the key ingredients of the group $\rG_n$ are
introduced, namely, the groups $S_n$, $\mT^n$ and $\Inn (\mI_n)$.
 The subgroup of $\rG_n$ they generate is their skew direct
 product $S_n\ltimes \mT^n\ltimes \Inn (\mI_n)$. An important
 description of the group $\rG_n$ is given (Corollary
 \ref{d21Mar9}.(2)), and two criteria of equality of two
 automorphisms of the algebra $\mI_n$ are obtained (Theorem
 \ref{21Mar9} and Corollary \ref{yz16Apr9}).

{\bf The group of inner automorphisms $\Inn (\mI_n)$ of the
algebra $\mI_n$}. For a group $G$, let $Z(G)$ denote its centre.
Since $\ga_n$ is an ideal of the algebra $\mI_n$, the intersection
$(1+\ga_n)^*:= \mI_n^*\cap (1+\ga_n)$ is a subgroup of the group
$\mI_n^*$ of units of the algebra $\mI_n$.

\begin{theorem}\label{A25Oct9}

\begin{enumerate}
\item {\rm (Theorem 4.5, \cite{algintdif})} $\mI_n^* = K^* \times
(1+\ga_n )^*$ and $Z(\mI_n^*) = K^*$. \item The map $(1+\ga_n )^*
\ra \Inn (\mI_n)$, $u\mapsto \o_u$, is a group isomorphism where
$\o_u(a) := uau^{-1}$, i.e. $\Inn (\mI_n)=\{ \o_u\, | \, u\in
(1+\ga_n )^*\}$. \item $\mI_1^* = K^* \times (1+F)^* \simeq K^*
\times \GL_\infty (K)$.
\end{enumerate}
\end{theorem}

{\it Proof}. Statement 2 follows from statement 1. Statement 3
follows from statement 1 and the fact that $(1+F)^* \simeq
\GL_\infty (K)$.  $\Box $

$\noindent $

{\bf The torus $\mT^n$}. The $n$-dimensional torus $\mT^n:= \{
t_\l \, | \, \l = (\l_1, \ldots , \l_n) \in K^{*n}\}$ is a
subgroup of the group $\rG_n$ where
$$t_\l (\int_i) = \l_i\int_i, \;\; t_\l (\der_i)= \l_i^{-1}\der_i, \;\; t_\l(H_i)= H_i, \;\; i=1,
\ldots , n.$$ Note that $t_\l (x_i) = \l_ix_i$ since $x_i=
\int_iH_i$. $\mT^n = \prod_{i=1}^n \mT^1(i)$ where $\mT^1(i) := \{
t_\l (i):= t_{(1, \ldots , 1, \l , 1, \ldots , 1)} \, | \, $ $\l
\in K^*\}$ and the scalar $\l$ is on $i$'th place.

$\noindent $

{\bf The symmetric group $S_n$}. The symmetric group $S_n$ is a
subgroup of the group $\rG_n$ where, for $\tau \in S_n$,
$$ \tau (\int_i) = \int_{\tau (i)}, \;\; \tau (\der_i) = \der_{\tau (i)}, \;\;
\tau (H_i) = H_{\tau (i)},\;\; i=1, \ldots , n.$$ The subgroup of
$\rG_n$ generated by $S_n$ and $\mT^n$ is the semi-direct product
$S_n\ltimes \mT^n$ since $S_n\cap \mT^n=\{ e\}$ and
\begin{equation}\label{tt1}
 \tau t_\l \tau^{-1}= t_{\tau (\l )}\;\; {\rm where}\;\; \tau (\l
) := (\l_{\tau^{-1}(1)},  \ldots , \l_{\tau^{-1}(n)}),
\end{equation}
 for
all $\tau \in S_n$ and $t_\l \in \mT^n$.

Since $\ga_n$ is the only maximal ideal of the algebra $\mI_n$,
$\s (\ga_n ) = \ga_n$ for all $\s \in \rG_n$. There is the group
homomorphism (recall that $B_n = \mI_n / \ga_n$):
\begin{equation}\label{Ixiaut}
\xi : \rG_n \ra \Aut_{K-{\rm alg}}(B_n), \;\; \s\mapsto
(\overline{\s} : a+\ga_n \mapsto \s (a) +\ga_n).
\end{equation}
The subgroup $S_n \ltimes \mT^n$ of $\rG_n$ maps isomorphically to
its image and $\xi (\Inn (\mI_n))= \{ e\}$, by Theorem
\ref{A25Oct9}.(2). Therefore, the subgroup $\rG_n'$ of $\rG_n$
generated by the subgroups $S_n$, $\mT^n$, and $\Inn (\mI_n)$ is
equal to their skew direct product, 
\begin{equation}\label{Gnpsk}
\rG_n'= S_n \ltimes \mT^n\ltimes \Inn (\mI_n).
\end{equation}
The goal of the paper is to prove that $\rG_n = \rG_n'$ (Theorem
\ref{25Oct9}).

$\noindent $

{\bf A description of the group $\rG_n$}. Let $A$ be an algebra
and $\s$ be its automorphism. For an $A$-module $M$, the {\em
twisted} $A$-module  ${}^{\s}M$, as a vector space,  coincides
with the module $M$ but the action of the algebra $A$ is given by
the rule: $a\cdot m :=\s (a)m$ where $a\in A$ and $m\in M$.
  The next lemma is
useful in finding the group of automorphisms of algebras that have
a {\em unique faithful} module satisfying  some
isomorphism-invariant properties.

\begin{lemma}\label{c21Mar9}
Suppose that an algebra $A$ has a unique (up to isomorphism)
faithful $A$-module $M$ that satisfies an isomorphism-invariant
property, say $\CP$. Then
$$ \Aut_{K-{\rm alg}}(A) = \{ \s_\v \, | \, \v \in \Aut_K(M), \;
\v A\v^{-1} = A\}$$ where  $\s_\v (a) := \v a \v^{-1}$ for $a\in
A$,  and the algebra $A$ is identified with its isomorphic copy in
$\End_K(M)$ via the algebra monomorphism $a\mapsto (m\mapsto am)$.
\end{lemma}

{\it Proof}.  Let $\s \in  \Aut_{K-{\rm alg}}(A)$. The twisted
$A$-module  ${}^{\s}M$ is faithful and satisfies the property
$\CP$. By the uniqueness of $M$, the $A$-modules $M$ and
${}^{\s}M$ are isomorphic. So, there exists an element $\v \in
\Aut_K(M)$ such that $\v a = \s (a)\v$ for all $a\in A$, and so
$\s (a) = \v a \v^{-1}$, as required. $\Box $

$\noindent $

{\it Example}. The matrix algebra $M_d(K)$ has a unique (up to
isomorphism)  simple module  which is $K^n$. Then, by Lemma
\ref{c21Mar9}, every automorphism of $M_d(K)$ is inner.

$\noindent $

Recall that the polynomial algebra $P_n$ is a unique (up to
isomorphism) faithful,  simple module for the algebra $\mI_n$
(Proposition \ref{c17Oct9}.(1)) and the algebra $\mA_n$,
\cite{jacaut}.
\begin{corollary}\label{d21Mar9}
\begin{enumerate}
\item {\rm \cite{jacaut}} $\mG_n = \{ \s_\v \, | \, \v \in
\Aut_K(P_n), \; \v \mA_n \v^{-1} = \mA_n\}$ where $\s_\v (a) := \v
a \v^{-1}$, $a\in \mA_n$. \item
 $\rG_n = \{
\s_\v \, | \, \v \in \Aut_K(P_n), \; \v \mI_n \v^{-1} = \mI_n\}$
where $\s_\v (a) := \v a \v^{-1}$, $a\in \mI_n$.
\end{enumerate}
\end{corollary}

In \cite{jacaut}, Corollary \ref{d21Mar9}.(1) was used in finding
the group $\mG_n$.

$\noindent $

{\bf The automorphism $\widehat{*}\in \Aut (\rG_n)$}. The
involution $*$ of the algebra $\mI_n$ induces the automorphism
$\widehat{*}$ of the group $\rG_n$ by the rule 
\begin{equation}\label{*GnI}
\widehat{*}:\rG_n\ra \rG_n, \;\; \s \mapsto *\circ \s \circ
*^{-1}.
\end{equation}

{\bf The $\mI_n$-module $P_n$}. By the very definition of the
algebra $\mI_n$ as a subalgebra of $\End_K(P_n)$, the
$\mI_n$-module $P_n$ is faithful. For the Weyl algebra $A_n$, the
$A_n$-module $A_n/\sum_{i=1}^n A_n\der_i$ is isomorphic to $P_n$
via $1+\sum_{i=1}^n A_n\der_i\mapsto 1$. The same statement is
true for the algebra $\mI_n$ (Proposition \ref{c17Oct9}.(3)).

\begin{proposition}\label{c17Oct9}
{\rm \cite{algintdif}}
\begin{enumerate}
\item  The polynomial algebra $P_n$ is the only (up to
isomorphism)  faithful, simple $\mI_n$-module. \item
$\mI_1=\mI_1\der \bigoplus \mI_1e_{00}$ and $\mI_1= \int
\mI_1\bigoplus e_{00}\mI_1$. \item ${}_{\mI_n}P_n\simeq \mI_n /
\sum_{i=1}^n \mI_n \der_i$.
\end{enumerate}
\end{proposition}
 The $\mI_n$-module $P_n$ is a very special module for the algebra
 $\mI_n$. Its properties, especially the uniqueness, are used
 often in this paper. The polynomial algebra $P_n=
 \bigoplus_{\alpha \in \N^n}Kx^\alpha$ is a naturally $\N^n$-graded
 algebra. This grading is compatible with the $\Z^n$-grading of
 the algebra $\mI_n$, i.e. the polynomial algebra $P_n$ is a $\Z^n$-graded $\mI_n$-module.
  Each element $\der_i\in \mI_n\subseteq
 \End_K(P_n)$ is a {\em locally nilpotent} map, that is
 $P_n=\bigcup_{j\geq 1}\ker_{P_n}(\der_i^j)$. Moreover,
 $$ \bigcap_{i=1}^n \ker_{P_n}(\der_i) = K.$$
Each element $\int_i\in \mI_n\subseteq
 \End_K(P_n)$ is an injective (but not a surjective) map. Each
 element $H_i\in \mI_n\subseteq
 \End_K(P_n)$ is a {\em semi-simple} map (that is $P_n=
 \bigoplus_{\l \in K} \ker_{P_n} (H_i-\l )$) with the set of
 eigenvalues $\Z_+:=\{ 1,2, \ldots \}$ since $H_i*x^\alpha =
 (\alpha_i+1) x^\alpha$ for all $\alpha \in \N^n$. Moreover,
\begin{equation}\label{kerHia}
\bigcap_{i=1}^n \ker_{P_n}(H_i-(\alpha_i+1))=Kx^\alpha,
 \;\;
 \alpha \in \N^n.
\end{equation}
In particular, the $K[H_1, \ldots , H_n]$-module
$P_n=\bigoplus_{\alpha \in \N^n}Kx^\alpha$ is the sum of simple,
non-isomorphic,  one-dimensional submodules $Kx^\alpha$, and so
$P_n$ is a semi-simple $K[H_1, \ldots , H_n]$-module.

\begin{corollary}\label{y16Apr9}
\begin{enumerate}
\item Let $M$ be an $\mI_n$-module. Then $\Hom_{\mI_n}(P_n,
M)\simeq \bigcap_{i=1}^n \ker (\der_{i,M})$, $f\mapsto f(1)$,
where $\der_{i,M}:M\ra M$, $m\mapsto \der_i m$. In particular,
$\End_{\mI_n}(P_n) \simeq K$.  \item By Proposition
\ref{c17Oct9}.(1,3), for each automorphism $\s \in \rG_n$, the
$\mI_n$-modules $P_n$ and ${}^\s P_n$ are isomorphic, and each
isomorphism $f:P_n\ra {}^\s P_n$ is
 given by the rule: $f(p) = \s (p) *v$, where  $v=f(1)$ is any
 nonzero element
 of the 1-dimensional vector space $\bigcap_{i=1}^n \ker (\s
 (\der_i)_{P_n})$.
\end{enumerate}
\end{corollary}

As an application of these results to the $\mI_n$-module $P_n$, we
have two useful  criteria of equality of two elements in the group
$\rG_n$.  They are   used in many proofs in this paper.

For an algebra $A$ and a subset $S\subseteq A$, $\Cen_A(S):=\{
a\in A\, | \, as=sa$ for all $s\in S\}$ is the {\em centralizer}
of the set $S$ in $A$. It is a subalgebra of $A$. It follows from
the presentation of the algebra $\mI_n$ as a GWA that
\begin{equation}\label{ACenxy}
\Cen_{\mI_n}(x_1, \ldots , x_n)=P_n, \;\; \Cen_{\mI_n}(\der_1,
\ldots , \der_n)=K[\der_1, \ldots , \der_n], \;\;
\Cen_{\mI_n}(\int_1, \ldots , \int_n)=K[\int_1, \ldots , \int_n].
\end{equation}
In more detail, since $\mI_n=\bigotimes_{i=1}^n\mI_1(i)$ it
suffices to prove the equalities for $n=1$, but in this case the
equalities are obvious.

Let $\mE_n:= \End_{K-{\rm alg}}(\mI_n)$ be the monoid of all the
$K$-algebra endomorphisms of $\mI_n$. The group of units of this
monoid is $\rG_n$. The automorphism $\widehat{*} \in \Aut (\rG_n)$
can be extended to an automorphism $\widehat{*} \in \Aut (\mE_n)$
of the monoid $\mE_n$: 
\begin{equation}\label{A1etah}
\widehat{*} : \mE_n\ra \mE_n, \;\; \s \mapsto *\circ \s \circ
*^{-1}.
\end{equation}
For each  element $\alpha \in \N^n$, let $x^{[\alpha
]}:=\int^\alpha *1$. Then   $x^{[\alpha ]}:=
\frac{x^\alpha}{\alpha !}:= \prod_{i=1}^n
\frac{x_i^{\alpha_i}}{\alpha_i !}$ and the set $\{ x^{[\alpha ]}\,
| \, \alpha \in \N^n\}$ is a $K$-basis for the polynomial algebra
$P_n$. The next result is instrumental in finding the group of
automorphisms of the algebra $\mI_n$.
\begin{theorem}\label{21Mar9}
{\rm (Rigidity  of the group $\rG_n$)} Let $\s , \tau \in \rG_n$.
Then the following statements are equivalent.
\begin{enumerate}
 \item $\s = \tau$. \item $\s (\int_1) = \tau (\int_1), \ldots , \s (\int_n) = \tau
 (\int_n)$.
 \item $\s (\der_1) = \tau (\der_1), \ldots , \s (\der_n) = \tau (\der_n)$.
 \item $\s (x_1) = \tau (x_1), \ldots , \s (x_n) = \tau (x_n)$.
\end{enumerate}
\end{theorem}

{\it Remark}.  It is not true that $\s (H_i) = \tau (H_i)$ for all
$i=1, \ldots , n$ implies $\s = \tau$ (Corollary
\ref{d9Oct9}.(2)).

 {\it Proof}. Without loss of generality we may assume that $\tau =e$, the
identity automorphism. The proof consists of two parts:
$(1\Leftrightarrow 2 \Leftrightarrow 3)$ and $(4\Rightarrow 1)$.
 Consider the following two subgroups of $\rG_n$, the stabilizers
of the sets $\{ \int_1, \ldots , \int_n\}$ and $\{ \der_1, \ldots
, \der_n\}$:
\begin{eqnarray*}
\St_{\rG_n}  (\int_1, \ldots , \int_n):= & \{ g\in \rG_n\, | \,
g(\int_1) = \int_1,
\ldots , g(\int_n)=\int_n\}, \\
 \St_{\rG_n}    (\der_1, \ldots , \der_n):= & \{ g\in \rG_n\,
| \, g(\der_1) = \der_1, \ldots , g(\der_n)=\der_n\}.
\end{eqnarray*}
Then $$\widehat{*} (\St_{\rG_n}    (\int_1, \ldots ,
\int_n))=\St_{\rG_n} (\der_1, \ldots , \der_n), \;\; \widehat{*}
(\St_{\rG_n} (\der_1, \ldots , \der_n))= \St_{\rG_n} (\int_1,
\ldots , \int_n).$$ Therefore, to prove that $(1\Leftrightarrow 2
\Leftrightarrow 3)$ (where $\tau = e$) is equivalent to show that
$\St_{\rG_n} (\int_1, \ldots , \int_n)=\{ e\}$. So, let $\s \in
\St_{\rG_n} (\int_1, \ldots , \int_n)$. We have to show that $\s =
e$, i.e. $\s (\der_i) = \der_i$ and $\s (H_i) = H_i$ for all $i$.
For each $i=1, \ldots , n$, $1=\s (\der_i\int_i) = \s (\der_i)
\int_i$ and $1=\der_i\int_i$. By taking the difference of these
equalities we see that $a_i:= \s (\der_i) -\der_i\in
\lann_{\mI_n}(\int_i)$. By (\ref{laI}), $a_i=\sum_{j\geq 0}
\l_{ij}e_{j0}(i)$ for some elements $\l_{ij}\in \bigotimes_{k\neq
i}\mI_1(k)$, and so $$ \s (\der_i) = \der_i+\sum_{j\geq 0}\l_{ij}
e_{j0}(i).$$ The element $\s (\der_i)$ commutes with the elements
$\s (\int_k) = \int_k$, $k\neq i$, hence all $\l_{ij}\in K[\int_1,
\ldots , \widehat{\int_i}, \ldots , \int_n]$, by (\ref{ACenxy}).
Since $e_{j0}(i) = \int_i^j e_{00}(i)$, we can write
$$ \s (\der_i) =\der_i+ p_ie_{00}(i) \;\; {\rm for \; some}\;\; p_i\in
K[\int_1, \ldots , \int_n].$$ We have to show that all $p_i=0$.
Suppose that this is not the case. Then $p_i\neq 0$ for some $i$.
We seek a contradiction. Note that $\s^{-1}\in \St_{\rG_n}
(\int_1, \ldots , \int_n)$, and so $\s^{-1} (\der_i) =\der_i+
q_ie_{00}(i)$ for some  $q_i\in K[\int_1, \ldots , \int_n]$.
Recall that $e_{00}(i) = 1-\int_i\der_i$. Then $\s^{-1}
(e_{00}(i))= 1-\int_i(\der_i+q_ie_{00}(i))=
(1-\int_iq_i)e_{00}(i)$, and
$$ \der_i=\s^{-1}\s (\der_i) = \s^{-1} (\der_i+p_ie_{00}(i))=
\der_i+(q_i+p_i(1-\int_iq_i))e_{00}(i), $$ and so $q_i+p_i=
\int_ip_iq_i$ since the map $K[\int_1, \ldots , \int_n]\ra
K[\int_1, \ldots , \int_n]e_{00}$, $p\mapsto pe_{00}$,  is a
bijection, by (\ref{Ieij}). This is impossible by comparing the
total degrees (with respect to the integrations) of the elements
on both sides of the equality. Therefore, $\s (\der_i) = \der_i$
for all $i$.

By Corollary \ref{y16Apr9}, there is an $\mI_n$-module isomorphism
$\v : P_n\ra {}^\s P_n$, $p\mapsto \s (p)*v=p*v=pv$, where $v:=\v
(1) \in \bigcap_{i=1}^n \ker_{{}^\s P_n}(\s (\der_i))=
\bigcap_{i=1}^n \ker_{P_n}(\der_i) = K1$.  Without loss of
generality we may assume that $v=1$. Then  $1=\v (1) = \v (H_i*1)
= \s (H_i) *1$ for all $i$.  For $i=1, \ldots , n$,
\begin{eqnarray*}
 \s (H_i) *x^{[\alpha ]} &=&\s (H_i) \int^\alpha   *1= \s (H_i) \s (\int^\alpha  ) *1= \s (H_i\int^\alpha )*1= \s(\int^\alpha  (H_i+\alpha_i))*1 \\
 &=& \s (\int^\alpha ) ( \s (H_i)+\alpha_i)*1=  \int^\alpha(\alpha_i+1)*1=
 (\alpha_i+1) x^{[\alpha ]}.
\end{eqnarray*}
This means that the linear maps $\s (H_i), H_i\in \End_K(P_n)$
coincide.  Therefore, $\s (H_i) = H_i$ for all $i$ since the
$\mI_n$-module $P_n$ is faithful. This proves that $\s =e$.

$(4\Rightarrow 1)$ Suppose that $\s (x_i) = x_i$ for all $i$. Then
$\s (p) = p$ for all polynomials $p\in P_n$. We have to show that
$\s = e$. By Corollary \ref{y16Apr9}, there exists the
$\mI_n$-module isomorphism $f: P_n\ra {}^\s P_n$, $p\mapsto \s
(p)*p = pv$. The map $f$ is a bijection hence $v\in K^*$. Without
less of generality we may assume that $v=1$, then $f = \id_{P_n}$.
Let $a\in \{ \der_i , H_i, \int_i\}$ and $b\in P_n$. Then
$$ a*b = f(a*b) = \s (a) *f(b) = \s (a) *b,$$
and so $\s (a) = a$ since the $\mI_n$-module $P_n$ is faithful.
This means that $\s =e$. The proof of the theorem is complete.
$\Box $

$\noindent $

Theorem \ref{21Mar9} means that $\St_{\rG_n}  (\int_1, \ldots ,
\int_n)= \St_{\rG_n}    (\der_1, \ldots , \der_n)= \St_{\rG_n}
(x_1, \ldots , x_n)=\{ e\}$.

$\noindent $

In zero characteristic, the Weyl algebra is the ring $\CD (P_n)$
of differential operators on the polynomial algebra $P_n$. In
prime characteristic, the Weyl algebra $A_n$ and the algebra $\CD
(P_n)$ are distinct, and the algebra $\CD (P_n)$ is much more
complicated object than the Weyl algebra $A_n$. An analogue of
Theorem \ref{21Mar9} does not hold for the algebra $\CD (P_n)$ in
characteristic zero, but does hold in prime characteristic
(Theorem 1.1, \cite{autlaur}). Also, the Rigidity Theorem is true
for the Jacobian algebra $\mA_n$ (Theorem 4.12, \cite{jacaut}) and
for the algebra $\mS_n$
 of one-sided inverses of the polynomial algebra $P_n$ (Theorem
 3.7, \cite{shrekaut}) but the Rigidity Theorem fails for the
 polynomial algebra $P_n$.

The ideal $F_n$ is the smallest nonzero ideal of the algebra
$\mI_n$. Therefore, $\s (F_n) = F_n$ for all $\s \in \rG_n$. The
next corollary shows that the action of the group $\rG_n$ on the
ideal $F_n$ is faithful. This result is used in the proof of the
fact that the group $\rG_n$ has trivial centre (Theorem
\ref{14Nov9}).

\begin{corollary}\label{yz16Apr9}
Let $\s , \tau \in \rG_n$. Then $\s = \tau$ iff $\s (e_{\alpha
0})=\tau (e_{\alpha 0})$ for all $\alpha \in \N^n$ iff $\s
(e_{0\alpha})=\tau (e_{0\alpha})$ for all $\alpha \in \N^n$ iff
$\s (e_{\alpha \beta })=\tau (e_{\alpha \beta})$ for all $\alpha ,
\beta \in \N^n$.
\end{corollary}

{\it Proof}.  The last `iff' follows from the previous two. The
second `iff' follows from the first one by using the automorphism
$\widehat{*}$ of the group $\rG_n$: $\s =\tau$ iff $\widehat{*}
(\s ) = \widehat{*} (\tau )$ iff $\widehat{*} (\s ) (e_{\alpha
0})= \widehat{*} (\tau ) (e_{\alpha 0})$ for all $\alpha \in \N^n$
(by the first `iff') iff $
 \s (e_{0\alpha})^*=\tau (e_{0\alpha})^*$ for all $\alpha \in
\N^n$ (since $e_{\alpha 0}^* = e_{0\alpha}$) iff $\s (e_{0\alpha
})= \tau (e_{0\alpha })$ for all $\alpha \in \N^n$.

So, it remains to prove that if $\s (e_{\alpha 0}) = \tau
(e_{\alpha 0})$ for all $\alpha \in \N^n$  then $\s = \tau$.
Without loss of generality we may assume that $\tau =e$, the
identity of the group $\rG_n$. So, we have to prove that if $\s
(e_{\alpha 0})= e_{\alpha 0}$ for all $\alpha \in \N^n$ then $\s
=e$. For each number $i=1, \ldots , n$,
$$ 0=(1-e_{00}(i))*1=\s (1-e_{00}(i))*1=\s (\int_i\der_i) *1= \s (\int_i)
\s(\der_i)*1, $$ and so $0=\s (\der_i) \s (\int_i) \s (\der_i) *1=
\s (\der_i\int_i) \s (\der_i) *1=\s (\der_i) *1$, i.e.
$\bigcap_{i=1}^n \ker (\s (\der_i)_{P_n})=K$. By Corollary
\ref{y16Apr9}.(2), the map $f:P_n\ra {}^\s P_n$, $p\mapsto \s (p)
*1$, is an $\mI_n$-module isomorphism. Now, $f(x^\alpha ) =
f(\alpha !e_{\alpha 0}*1)= \s (\alpha !e_{\alpha 0})*1=\alpha
!e_{\alpha 0}*1=x^\alpha$ for all $\alpha \in \N^n$ where $\alpha
!:= \alpha_1 !\cdots \alpha_n !$. This means that $f$ is the
identity map. For all $a\in \mI_n$ and $p\in P_n$, $a*p= f(a*p)=
\s (a)*f(p)= \s (a) *p$, and so $\s (a) = a$ since the
$\mI_n$-module $P_n$ is faithful. That is $\s = e$, as required.
$\Box $
\begin{corollary}\label{a29Oct9}

\begin{enumerate}
\item  Let $\ga$ be a nonzero ideal of the algebra $\mI_n$ and $\s
, \tau \in \rG_n$. Then $\s = \tau$ iff $\s (a) = \tau (a)$ for
all $a\in \ga$. \item Let $\ga$ be a nonzero ideal of the algebra
$\mA_n$ and $\s , \tau \in \mG_n$. Then $\s = \tau$ iff $\s (a) =
\tau (a)$ for all $a\in \ga$.
\end{enumerate}
\end{corollary}

{\it Proof}. 1. Since $F=\bigoplus_{\alpha, \beta \in \N^n}
Ke_{\alpha \beta}\subseteq \ga$, statement 1 follows from
Corollary \ref{yz16Apr9}.

2. Similarly, since $F=\bigoplus_{\alpha, \beta \in \N^n}
Ke_{\alpha \beta}\subseteq \ga$, statement 2 follows from
(Corollary 4.6, \cite{jacaut}). $\Box$


\section{The group $\Aut_{K-{\rm alg}}(\mI_1)$}\label{AAII1}

In this section, the group $\rG_1$ and its explicit generators are
found (Theorem \ref{31Oct9}).  The key idea in finding the group
$\rG_1$ of automorphisms of the algebra $\mI_1$ is to use Theorem
\ref{21Mar9}, some of the properties of the index of linear maps
in the vector space $P_1=K[x]$, and the
 explicit structure  of the group $\Aut_{K-{\rm alg}}(\mI_1/ F)$ (Theorem
\ref{A31Oct9}). It is proved that any algebra endomorphism of the
algebra $\mI_1$ is a monomorphism (Theorem \ref{B31Oct9}); note
that the algebra $\mI_1$ is not a simple algebra.

$\noindent $

{\bf The group $\brG_1:= \Aut_{K-{\rm alg}}(B_1)$}. Recall that
$B_1= K[H][x,x^{-1}; \s ]$ and $\s (H) = H-1$. Consider the
following automorphisms of the algebra $B_1$:
\begin{eqnarray*}
 t_\l:& x\mapsto \l x, &\;\; H\mapsto H, \;\;\;\;\;\;\;\;  (\l \in K^*) \\
 s_p:&   x\mapsto x, &\;\; H\mapsto H+p, \;\;\; (p\in K[x,x^{-1}])  \\
 \zeta :&  x\mapsto x^{-1},& \;\; H\mapsto - H,
\end{eqnarray*}
and the subgroups they generate in the group $\brG_1$:
$$\mT^1:=  \{ t_\l \, | \, \l\in K^*\}\simeq K^*,\;\;  \Sh_1:= \{ s_p \, | \, p\in
K[x,x^{-1}]\} \simeq K[x,x^{-1}], \;\; \langle \zeta \rangle
\simeq \Z_2:= \Z / 2\Z .$$
 We can easily check that
\begin{equation}\label{2Hpm0}
\zeta t_\l\zeta^{-1} = t_\l^{-1},  \;\;  \zeta s_p\zeta^{-1} =
s_{-\zeta (p)}, \;\; t_\l s_pt_\l^{-1} = s_{t_\l (p)}.
\end{equation}
 It follows that the subgroup of $\brG_1$  generated by  the three
 subgroups above is, in fact, their semidirect product,
$$ \langle \zeta \rangle \ltimes \mT^1\ltimes
\Sh_1\simeq  \Z_2\ltimes  K^*\ltimes K[x,x^{-1}],$$ since, for
$\varepsilon =0,1$; $\l \in K^*$;
  and $p\in K[x,x^{-1}]$:
\begin{equation}\label{4sthp}
 \zeta^\varepsilon t_\l s_p : x\mapsto \l x^{1-2\varepsilon},
\;\; H\mapsto (-1)^\varepsilon H+\zeta^\varepsilon t_\l (p),
\end{equation}
 and  $\zeta^\varepsilon t_\l s_p= e$ iff $\varepsilon =0$, $\l =0$,
 and $p=0$. Theorem \ref{A31Oct9} shows that this is the whole
 group of automorphisms of the algebra $B_1$.

For a group $G$, $[G,G]$ denotes its {\em
 commutant}, i.e. the subgroup of $G$ generated by all the group {\em
 commutators} $[a,b]:=aba^{-1}b^{-1}$ of the elements $a,b\in G$.
 The centre of a group $G$ is denoted by $Z(G)$. For subgroups $A$
 and $B$ of the group $G$, let $[A, B]$ be the subgroup of $G$
 generated by all the commutators $[a,b]$ where $a\in A$ and $b\in
 B$. Given a skew direct product $A\ltimes \prod_{i=1}^m B_i$ of
 groups such that $aB_ia^{-1}\subseteq B_i$ for all $a\in A$ and
 $i=1, \ldots , m$; then its commutant is equal to $[A,A]\ltimes
 \prod_{i=1}^m ([A, B_i]\cdot [B_i , B_i])$ (Lemma 5.4.(2),
 \cite{jacaut}). This fact is used in the proof of the following
 theorem.

\begin{theorem}\label{A31Oct9}
\begin{enumerate}
\item $\brG_1=\langle \zeta \rangle \ltimes \mT^1\ltimes \Sh_1$.
\item $Z(\brG_1) = \{ e\}$. \item $[\brG_1,\brG_1]=\{ t_{\l^2}\, |
\, \l \in K^*\}\ltimes \Sh_1$ and $\brG_1/[\brG_1,\brG_1]\simeq
\Z_2\times K^*/K^{*2}$.
\end{enumerate}
\end{theorem}

{\it Proof}. 1. Let $\s \in \brG_1$ and $\brG_1'$ be the skew
direct product. It suffices to show that $\s \in \brG_1'$ since
$\brG_1'\subseteq \brG_1$. The automorphism $\s $ of the algebra
$B_1$ induces an automorphism of its group of units $B_1^*
=\bigcup_{i\in \Z}K^*x^i$. Then $\s (x) =\l x^{\pm 1}$ for a
nonzero scalar $\l \in K^*$. Multiplying $\s $ on the left by a
suitable element of the group $\langle \zeta \rangle \ltimes
\mT^1$ we may assume that $\s (x) = x$. Then
$$[\s (H)-H, x]= \s ([H,x]) -[H,x]=\s (x) - x=0.$$ Therefore,
$p:=\s (H) - H\in \Cen_{B_1}(x)=K[x,x^{-1}]$, and so $\s = h_p\in
\brG_1'$. This proves that $\brG_1 = \brG_1'$.

2. Let $z\in Z (\brG_1)$. By statement 1, $z=\zeta^\varepsilon
t_\l s_p$ for some elements $\varepsilon = 0,1$; $\l \in K^*$; and
$p\in K[x,x^{-1}]$.  By (\ref{2Hpm0}),  $\zeta^\varepsilon t_{\l
\mu}s_{t_{\mu^{-1}}(p)}= zt_\mu = t_\mu z=\zeta^\varepsilon t_{\l
\mu^{1-2\varepsilon}}s_p$, hence $\varepsilon =0$ and $p\in K$.
Next, $\zeta t_\mu s_p = \zeta z = z\zeta = \zeta
t_{\mu^{-1}}s_{-p}$, hence $\mu = \pm 1$ and $p=0$, i.e. $z=t_{\pm
1}$. Since $s_xt_{-1}\not=t_{-1}s_x$, $z=t_1=e$. Therefore,
$Z(\brG_1) = \{ e\}$.

3. It suffices to prove only that the equality holds since then
the isomorphism is obvious, by statement 1. Let $R$ be the RHS of
the equality. Then $R\subseteq [ \brG_1, \brG_1]$ since
$$ t_{\l^2} = [ \zeta , t_{\l^{-1}}], \;\; s_{\mu x^j} = [ t_2,
s_{\frac{\mu x^j}{2^j-1}}], \;\; s_\mu = [ \zeta ,
s_{-\frac{\mu}{2}}], $$ where  $0\neq j\in \Z$, $\l \in K^*$, and
$\mu \in K$. It suffices to show that $[ \brG_1/\Sh_1,
\brG_1/\Sh_1]\subseteq R'$ where $R':= \{ t_{\l^2}\, | \, \l \in
K^*\}$ is treated as a subgroup of the factor group
$\brG_1/\Sh_1\simeq \langle \zeta \rangle\ltimes \mT^1$. By Lemma
5.4, \cite{jacaut}, $[\langle \zeta \rangle \ltimes \mT^1, \langle
\zeta \rangle \ltimes \mT^1]=[\langle \zeta \rangle, \mT^1]=R'$.
$\Box $

$\noindent $

{\bf The index $\ind$ of linear maps and Fredholm operators}. Let
$\CC = \CC (K)$ be the family of all $K$-linear maps with finite
dimensional kernel and cokernel (such maps are called the {\em
Fredholm linear maps/operatorts}).  So, $\CC$ is the family of
{\em Fredholm} linear maps/operators. For vector spaces $V$ and
$U$, let $\CC (V,U)$ be the set of all the linear maps from $V$ to
$U$ with finite dimensional kernel and cokernel. So, $\CC
=\bigcup_{V,U}\CC (V,U)$ is the disjoint union.

$\noindent $

{\it Definition}. For a linear map $\v \in \CC$, the integer
$$ \ind (\v ) := \dim \, \ker (\v ) - \dim \, \coker (\v )$$
is called the {\em index} of the map $\v$.

$\noindent $

{\it Example}. Note that $\der, \int \in \mI_1\subset
\End_K(P_1)$. Then  
\begin{equation}\label{indId}
\ind (\der^i)= i\;\; {\rm and }\;\; \ind (\int^i)= -i, \;\; i\geq
1.
\end{equation}

Each nonzero element $u$ of the skew Laurent polynomial algebra
$\CA_1= \CL_1[x,x^{-1}; \s_1 ]$  (where $\s_1(H)=H-1$) is a unique
sum $u=\l_sx^s+\l_{s+1}x^{s+1}+\cdots +\l_dx^d$ where all $\l_i\in
\CL_1$, $\l_d\neq 0$, and $\l_dx^d$ is the {\em leading term} of
the element $u$. Recall that $\CL_1:= K[H^{\pm 1}, (H\pm 1)^{-1},
(H\pm  2)^{-1}, \ldots ]$, $B_1\subset \CA_1$, and  $\mI_1\subset
\mA_1$. The integer $\deg_x(u)=d$ is called the {\em degree} of
the element $u$, $\deg_x(0):=-\infty $. For all $u,v\in \CA_1$,
$\deg_x(uv)= \deg_x(u) + \deg_x(v)$ and $\deg_x(u+v) \leq \max \{
\deg_x(u), \deg_x(v)\}$. The next lemma explains how to compute
the index of the elements $\mA_1\backslash F$ (resp.
$\mI_1\backslash F$) via the degree function $\deg_x$ and proves
that the index is a $\mG_1$-invariant (resp. a $\rG_1$-invariant)
concept. Note that $F\cap \CC = \emptyset$.
\begin{lemma}\label{c24Mar9}
\begin{enumerate}
\item {\rm \cite{jacaut}} $\CC \cap \mA_1= \mA_1\backslash F $
(recall that $\mA_1\subset \End_K(P_1)$) and, for each element
$a\in \mA_1\backslash F$, $ \ind (a) = -\deg_x(\oa )$ where $\oa =
a+F\in \mA_1/F= \CA_1$. \item \cite{jacaut} $\ind (\s (a)) = \ind
(a)$ for all $\s \in \mG_1$ and $a\in \mA_1\backslash F$.\item
$\CC \cap \mI_1= \mI_1\backslash F $; for each element $a\in
\mI_1\backslash F$, $ \ind (a) = -\deg_x(\oa )$ where $\oa =
a+F\in B_1$; and $\ind (\s (a)) = \ind (a)$ for all $\s \in \rG_1$
and $a\in \mI_1\backslash F$.
\end{enumerate}
\end{lemma}

{\it Proof}. Since $\mI_n\subseteq \mA_n$, statement 3 follows
from statements 1, 2 and Corollary \ref{d21Mar9}.(2)). $\Box $

$\noindent $

The next theorem presents the group $\rG_1$ and its explicit
generators.

\begin{theorem}\label{31Oct9}
\begin{enumerate}
\item $\rG_1=\mT^1\ltimes \Inn (\mI_1)$. \item $\rG_1\simeq
K^*\ltimes \GL_\infty (K)$. \item $[\rG_1, \rG_1] = \{ \o_u\, | \,
u\in \SL_\infty (K)\}$ and $\rG_1/[\rG_1, \rG_1] \simeq
\mT^1\times \mT^1$. \item The group $\rG_1$ is generated by the
elements $t_\l $, $\o_{1+\l e_{ij}}$ where $i\neq j$ and $\l \in
K^*$, and $\o_{1+\mu e_{11}}$ where $ \mu \in K\backslash \{ -1\}
$.
\end{enumerate}
\end{theorem}

{\it Proof}. 1. Let $\s \in \rG_1$. By (\ref{Gnpsk}), $\rG_1'=
\mT^1 \ltimes \Inn (\mI_1)\subseteq \rG_1$. It remains to show
that the reverse inclusion holds, that is $\s \in \rG_1'$. The
ideal $F$ of the algebra $\mI_1$ is the only maximal ideal.
Therefore, $\s (F) = F$ and $\overline{\s} :=\xi (\s )\in \brG_1$,
see (\ref{Ixiaut}). By Theorem \ref{A31Oct9}.(1) and
(\ref{4sthp}), either $\overline{\s}(\overline{\der} ) = \l
x^{-1}$ or, otherwise, $\overline{\s}(\overline{\der} ) = \l x$
for some element $\l \in K^*$. Equivalently, either $\s (\der ) =
\l \der +f$ or $\s (\der ) = \l \int +f$ for some element $f\in
F$. Recall that $\mI_1\subset \mA_1$. By Lemma \ref{c24Mar9}.(3)
and Corollary \ref{d21Mar9}.(2), the second case is impossible as
we have the contradiction:
$$ 1=\ind (\der ) = \ind (\s (\der ))= \ind (\l \int +f) =
\ind (\l xH^{-1} +f)= -\deg_x(\l xH^{-1}) = -1.$$ Therefore, $\s
(\der ) = \l \der +f$. Replacing $\s$ by $t_\l \s$ we may assume
that $\s (\der ) =  \der +g$ where $g:=t_\l (f) \in F$ (as $F$ is
the only maximal ideal of the algebra $\mI_1$, hence $\tau (F)=F$
for all $\tau \in \rG_1$). Fix a natural number $m$ such that
$g\in \sum_{i,j=0}^mKe_{ij}$. Then the finite dimensional vector
spaces
$$ V:= \bigoplus_{i=0}^m Kx^{[i]}\subset V':=
\bigoplus_{i=0}^{m+1}Kx^{[i]}$$ are $\der'$-invariant where
$\der':= \s (\der )= \der +g$,  $x^{[i]}:=\frac{x^i}{i!}$, and
$x^{[0]}:=1$. Note that $\der'*x^{[m+1]} = \der *x^{[m+1]}=
x^{[m]}$ since $g*x^{[m+1]}=0$.
 Note that $P_1=\bigcup_{i\geq 1}\ker (\der^i)$ and $\dim \, \ker_{P_1} (\der )
 =1$. Since the $\mI_1$-modules $P_1$ and ${}^{\s}P_1$ are
 isomorphic (Corollary \ref{y16Apr9}.(2)), $P_1=\bigcup_{i\geq 1}\ker (\der'^i)$
 and $\dim \, \ker_{P_1} (\der')
 =1$. This implies that the elements $x'^{[0]}, x'^{[1]}, \ldots , x'^{[m]},
 x^{[m+1]}$ are a $K$-basis for the vector space $V'$ where
 $$ x'^{[i]}:= \der'^{m+1-i}*x^{[m+1]}, \;\; i=0,1, \ldots , m; $$
 and the elements  $x'^{[0]}, x'^{[1]}, \ldots , x'^{[m]}$ are a $K$-basis for the vector space
 $V$. Then the elements
 $$ x'^{[0]}, x'^{[1]}, \ldots , x'^{[m]}, x^{[m+1]}, x^{[m+2]},  \ldots  $$
are a $K$-basis for the vector space $P_1$. The $K$-linear map
\begin{equation}\label{chm}
\v : P_1\ra P_1, \;\; x^{[i]}\mapsto x'^{[i]}\; (i=0,1,\ldots ,
m), \; x^{[j]}\mapsto x^{[j]}\; (j>m),
\end{equation}
belongs to the group $(1+F)^*= \GL_\infty (K)\simeq \Inn (\mI_1)$
(by Theorem \ref{A25Oct9}.(2)) and satisfies the property that $
\der'\v  = \v \der  $, the equality is in $\End_K(P_1)$. This
equality can be rewritten as follows:
$$ \o_{\v^{-1}}\s (\der ) = \der \;\; {\rm where}\;\;
\o_{\v^{-1}}\in \Inn (\mI_1).$$ By Theorem \ref{21Mar9},
$\o_{\v^{-1}}\s =e$, hence $\s \in \rG_1'$.

2. Statement 2 follows from statement 1  and the fact that $\Inn
(\mI_1)\simeq (1+F)^*  \simeq \GL_\infty (K)$ (by
 Theorem \ref{A25Oct9}).

3. $[\rG_1, \rG_1]=
 [\mT^1\ltimes  \GL_\infty (K), \mT^1\ltimes  \GL_\infty (K)]=
[\mT^1 , \GL_\infty (K)] [ \GL_\infty (K), \GL_\infty
(K)]=\SL_\infty (K)$ since $[\mT^1 , \GL_\infty (K)]\subseteq
\SL_\infty (K)$ and  $\SL_\infty (K) = [ \GL_\infty (K),
\GL_\infty (K)]$. Now, $\rG_1/ [\rG_1, \rG_1]\simeq \mT^1\times
\GL_\infty (K) / \SL_\infty (K) \simeq \mT^1 \times \mT^1$.

4. Statement 4 follows from statements 1 and  2 and  the fact that
the group $\GL_\infty (K)$ is generated by the elements $1+\l
e_{ij}$ and $1+\mu e_{11}$ where $i\neq j$, $\l \in K^*$, and $
\mu \in K\backslash \{ -1\}$. $\Box $

\begin{corollary}\label{a31Oct9}
$\xi (\rG_1) = \mT^1$ and $\ker (\xi ) = \Inn (\mI_1)$.
\end{corollary}

{\it Proof}. The homomorphism $\xi$   maps isomorphically the
torus $\mT^1$ onto its image, and $\xi (\Inn (\mI_1)) = \{ e\}$
(by Theorem \ref{A25Oct9}). Therefore, $\xi (\rG_1) = \mT^1$ and
$\ker (\xi ) = \Inn (\mI_1)$ since $\rG_1= \mT^1\ltimes \Inn
(\mI_1)$. $\Box $

$\noindent $

Every algebra endomorphism of a {\em simple} algebra is a {\em
monomorphism}. The algebra $\mI_1$ is not simple but the same
result holds.

\begin{theorem}\label{B31Oct9}
Every algebra endomorphism of the algebra $\mI_1$ is a
monomorphism.
\end{theorem}

{\it Proof}. Recall  that $F$ is the  only proper  ideal of the
algebra $\mI_1$, and $\mI_1/F =B_1:= K[H][x,x^{-1}; \s_1]$ is a
simple algebra where $\s_1(H) = H-1$. Suppose that $\g$ is an
algebra endomorphism of $\mI_1$ which is not a monomorphism, then
necessarily $\g (F)=0$, and the endomorphism $\g$ induces the
algebra monomorphism $\overline{\g}: B_1\ra \mI_1$, $a+F\mapsto \g
(a)$. We seek a contradiction. Since $\der \int=1$ and $\int \der
=1-e_{00}$, we have  the equalities $\g (\der) \g (\int ) =1$ and
$\g (\int ) \g (\der) =1$, i.e. the elements $\g (\der) $ and $\g
(\int )$ are units of the algebra $\mI_1$. Therefore, the images
of the elements  $\g (\der) $ and $\g (\int )$ in the algebra
$B_1$  under the epimorphism $\pi :\mI_1\ra B_1$ belong to the
group of units of the algebra $B_1$ which is  $K^*$, hence $\pi
(\im (\g ))\subseteq K\langle \pi \g (H)\rangle$, a commutative
algebra. This is impossible since the algebra $\im (\g )  \simeq
B_1$ is a simple non-commutative algebra. This contradiction
proves the theorem.
 $\Box $

$\noindent $

{\it Question. Is an algebra endomorphism of the algebra $\mI_1$
an isomorphism? The same question we can ask for $\mI_n$, see
Theorem \ref{B14Nov9}.}

$\noindent $

This question has flavour of the Question/Conjecture of Dixmier
\cite{Dix}: {\it is an algebra endomorphism of the Weyl algebra an
isomorphism?}


\section{The group of automorphisms
of the algebra $\mI_n$}\label{GAAIN}

In this section, it is proved that $\rG_n = S_n \ltimes \mT^n
\ltimes \Inn (\mI_n)$ (Theorem \ref{25Oct9}), the groups $\Inn
(\xi )$ and $\ker (\xi )$ are found (Theorem \ref{25Oct9}.(3,4)).
The group $\rG_n$ has trivial centre (Theorem \ref{14Nov9}).
 For each automorphism $\s$ of the
algebra $\mI_n$, an explicit inversion formula  is given via the
elements$\{ \s (\der_i), \s (\int_i) \, | \, i=1, \ldots , n\}$
(Theorem \ref{13Nov9}). It is proved that no proper prime factor
algebra of the algebra $\mI_n$ can be embedded into $\mI_n$
(Theorem \ref{B14Nov9}). It is shown that each automorphism of the
algebra $\CI_n$ of scalar integro-differential operators can be
uniquely extended to an automorphism of the algebra $\mI_n$
(Theorem \ref{2Nov9}).

$\noindent $

 {\bf The group
$\brG_n:= \Aut_{K-{\rm alg}}(B_n)$}. Recall that $B_n=K[H_1,
\ldots , H_n][z_1^{\pm 1}, \ldots , z_n^{\pm 1}; \s_1, \ldots ,
\s_n]$,  $\s_i (H_j) = H_j-\d_{ij}$. The algebra $B_n$ contains
the Laurent polynomial algebra $L_n:=K[z_1^{\pm 1}, \ldots ,
z_n^{\pm 1}]$. The set $L_n':= \{ (p_i)\in L_n^n\, | \,
z_i\frac{\der p_j}{\der z_i}=z_j\frac{\der p_i}{\der z_j}, \forall
i\neq j\}$ is a $K$-subspace of $L_n^n$. Clearly, $K^n\subseteq
L_n'$ and $L_n'\supseteq \bigoplus_{i=1}^n K[z_i, z_i^{-1}]=\{
(p_i)\, | \, p_i\in K[z_i, z_i^{-1}]\}$.  {\em The elements
$H_1+p_1, \ldots , H_n+p_n$ (where $p_i\in L_n$) of the algebra
$B_1$ commute iff} $(p_i) \in L_n'$.

Consider the following automorphisms of the algebra $B_n$: for
$i=1, \ldots , n$,
\begin{eqnarray*}
 a:& z_i\mapsto& \prod_{j=1}^n z_j^{a_{ij}}, \;\; H_i\mapsto \sum_{j=1}^n H_jb_{ji}, \;\;
  (a=(a_{ij})\in \GL_n(\Z ), (b_{ij})=a^{-1}) \\
 t_\l :& z_i\mapsto & \l_iz_i, \;\;\;\;\;\;\;\, H_i\mapsto H_i, \;\;\;\;\;\;\;\;\;\;\;\; (\l =
 (\l_i)\in K^{*n})\\
 s_p:& z_i\mapsto & z_i,\;\;\;\;\;\;\;\;\;  \;\; H_i\mapsto H_i+p_i, \;\;\;\;\; (p=(p_i)
 \in L_n')
\end{eqnarray*}
and the subgroups they generate in the group $\brG_n$:
$$\O_n:=\{ a\, | \, a\in \GL_n (\Z )\} \simeq \GL_n(\Z )^{op},
\;\; \mT^n := \{ t_\l \, | \, \l \in K^{*n}\}\simeq K^{*n}, \;\;
\Sh_n := \{ s_p\, | \, p\in L_n'\} \simeq L_n'.  $$  The group
$\O_n$  is isomorphic to the {\em opposite group} $\GL_n(\Z
)^{op}$ of the general linear group $\GL_n (\Z )$ via $a\mapsto
a$. Recall that  as a set $\GL_n(\Z )^{op}=\GL_n (\Z )$ but the
group structure on $\GL_n(\Z )^{op}$ is given by the rule $a\circ
b= ba$, the matrix multiplication. The group  $\GL_n(\Z )^{op}$ is
isomorphic to the group $\GL_n(\Z )$ via $a\mapsto a^{-1}$.

For each nonzero element $\alpha \in \Z^n$, the set $\Supp (\alpha
) :=\{ i\, | \, \alpha_i\neq 0\}$ is called the {\em support} of
$\alpha$, and $\min (\alpha )$ denotes the minimal number in the
support of $\alpha$. The following lemma gives a $K$-basis for the
vector space $L_n'$. Recall that $K^n \subseteq L_n'$.

\begin{lemma}\label{a8Nov9}
$L_n'=K^n\bigoplus \bigoplus_{0\neq \alpha \in \Z^n}Kb_\alpha$
where $b_\alpha=(\l_iz^\alpha)$, $\l_i =\begin{cases}
0& \text{if }i\not\in \Supp (\alpha ),\\
1& \text{if }i=\min (\alpha ), \\
\frac{\alpha_i}{\alpha_{\min (\alpha)}}& \text{if }i\in \Supp
(\alpha ).
\end{cases}$
\end{lemma}

{\it Proof}. It is obvious that $L_n'\supseteq R$ where $R$ is the
RHS of the equality. Each direct summand $Kz^\alpha$ of the
Laurent  polynomial algebra $K[z_1^{\pm 1}, \ldots , z_n^{\pm
1}]=\bigoplus_{\alpha\in \Z^n}Kz^\alpha$ is invariant under the
actions of the $K$-derivations $z_1\frac{\der}{\der z_1}, \ldots
,z_n\frac{\der}{\der z_n}$ since $z_i\frac{\der (z^\alpha)}{\der
z_i}=\alpha_iz^\alpha$. Therefore, a $K$-basis for the vector
space $L_n'$ can be chosen in such a way that every element of the
basis is of the type $(\l_iz^\alpha )$ where $0\neq (\l_i) \in
K^n$, $\alpha \in \Z^n$ and 
\begin{equation}\label{ailj}
\alpha_i\l_j=\alpha_j\l_i \;\; {\rm for \; all}\;\;  i\neq j.
\end{equation}
If $\alpha =0$ then there is no restriction on the scalars $\l_j$
and we get the vector space $K^n$.

If $\alpha \neq 0$ then either $\l_{\min (\alpha )}\neq 0$ or
$\l_{\min (\alpha )}= 0$. In the first case, the space of
solutions to the system of linear equations (\ref{ailj}) is
$K(\l_i)$ where the vector $(\l_i)\in K^n$ is as in the lemma. If
$\l_{\min (\alpha )}= 0$ then $\l_i=0$ for all $i=1, \ldots , n$.
Now, the lemma is obvious. $\Box $

$\noindent $

We can easily verify that 
\begin{equation}\label{atas}
at_\l a^{-1} =t_{\l^{a^{-1}}}, \;\; as_pa^{-1} = s_{a(p)a}, \;\;
t_\l s_p t_\l^{-1} = s_{t_\l (p)},
\end{equation}
where $\l^{a^{-1}}=(\l_i')$, $\l_i':=\prod_{j=1}^n \l_j^{b_{ij}}$,
$(b_{ij})= a^{-1}$; $a(p)a= (p_i')$, $p_i'=\sum_{j=1}^n
a(p_j)a_{ji}$; $t_\l (p) = (t_\l (p_i))$. It follows that the
subgroup of the group $\brG_n$ generated by the three subgroups
above is, in fact, the semi-direct product $\O_n \ltimes \mT^n
\ltimes \Sh_n$ since 
\begin{equation}\label{atsp}
at_\l s_p : z_i\mapsto \l_i\prod_{j=1}^n z_j^{a_{ij}}, \;\;
H_i\mapsto \sum_{j=1}^n H_jb_{ji}+at_\l (p_i),
\end{equation}
and so $at_\l s_p = e$ iff $a=e$, $\l =0$, and $p=0$. Theorem
\ref{1Nov9} shows that this is the whole group of  automorphisms
of the algebra $B_n$.

\begin{theorem}\label{1Nov9}
\begin{enumerate}
\item $\brG_n = \O_n \ltimes \mT^n \ltimes \Sh_n\simeq \GL_n(\Z
)^{op} \ltimes K^{*n}\ltimes L_n'$. \item $Z(\brG_n) = \{ e\}$.
\item Let $n\geq 2$. Then $[\brG_n , \brG_n] = [\O_n ,
\O_n]\ltimes \mT^n \ltimes \Sh_n$ and $\brG_n /[\brG_n ,
\brG_n]\simeq \Z_2$.
\end{enumerate}
\end{theorem}

{\it Proof}. 1. Let $\s \in \brG_n$ and $\brG_n'$ be the skew
direct product $\O_n \ltimes \mT^n \ltimes \Sh_n$. Recall that
$\brG_n'\subseteq \brG_n$. It remains to show that the reverse
inclusion holds. The automorphism $\s$ of the algebra $B_n$
induces an automorphism of its group of units
$B_n^*=\bigcup_{\alpha \in \Z^n} K^* z^\alpha$. Then $\s (z_i) =
\l_i\prod_{j=1}^n z_j^{a_{ij}}$ for $i=1, \ldots , n$, where
$\l_i\in K^*$ and $a=(a_{ij})\in \GL_n (\Z)$. Replacing the
automorphism $\s$ with $t_\mu a^{-1} \s$ for some $\mu \in K^{*n}$
we may assume that $\s (z_i) = z_i$ for all $i=1, \ldots , n$.
Then, for all indices $i,j=1, \ldots , n$,
$$ [ \s (H_i) - H_i, z_j]=\s ([H_i, z_j])-[H_i, z_j]= \s
(\d_{ij}z_j) - \d_{ij}z_j = \d_{ij} z_j - \d_{ij}z_j =0.$$
Therefore, $p_i := \s (H_i) - H_i \in \Cen_{B_n}(z_1, \ldots ,
z_n)= L_n$.  The elements $\s (H_1),  \ldots , \s (H_n)$ commute,
$$0=[\s (H_i), \s (H_j)]= [H_i+p_i, H_j+p_j]=[H_i, p_j]-[H_j,
p_i]=z_i\frac{\der p_j}{\der z_i}-z_j\frac{\der p_i}{\der z_j}.$$
  Therefore, $(p_i)\in L_n'$, i.e. $\s = s_p \in
\brG_n'$, and so $\brG_n \subseteq \brG_n'$.

2. By Theorem \ref{A31Oct9}.(2), we may assume that $n\geq 2$. Let
$z\in Z(\brG_n)$. By statement 1, $z=at_\l s_p$. By (\ref{atas}),
for all elements $b\in \O_n$, $bat_\l s_p = bz=zb=
abt_{\l^b}s_{b^{-1}(p)b^{-1}}$. Then $\l =(1, \ldots , 1)$  and
$b\in Z(\O_n) = \{ \pm e\}$ where $-e:z_i \mapsto z_i^{-1}$,
$H_i\mapsto -H_i$, for all $i=1, \ldots , n$. Since
$(-e)s_ps_{(z_1, \ldots , z_n)}\neq s_{(z_1, \ldots ,
z_n)}(-e)s_p$, we have $z=s_p$. The equalities $s_{t_\l (p)}= t_\l
s_p t^{-1}_\l = s_p$ for all $t_\l \in \mT^n$ imply that $p\in
K^n$. The equalities $s_{pa}=s_{a(p)a}= as_pa^{-1}= s_p$ for all
elements $a\in \O_n$ imply that $p=0$, and so $z=e$. Therefore,
$\Z (\O_n) = \{ e\}$.

3.  Let $R$ be the RHS of the first equality in statement 3. Then
$R \subseteq [ \brG_n , \brG_n ] $ since $[E-E_{ij}, t_\l
(j)]=t_\l (i)$ where $t_\l (i)\in \mT^1(i)$, $E\in M_n(K)$ is the
identity matrix, $E_{ij}\in M_n(K)$ are the matrix units, and
$$[t_2(i), s_{
\frac{b_\alpha}{2^{\alpha_i}-1} }]= s_{b_\alpha }, \;\;\;\; [
(-e), s_{-\frac{\mu }{2}e_j }]=s_{\mu e_j}$$ where $\alpha \in
\Z^n$ with $\alpha_i\neq 0$,  and the set $\{ e_1, \ldots , e_n\}$
is the standard  $K$-basis for the vector space $K^n$. The reverse
inclusion $R\supseteq [ \brG_n , \brG_n]$ is obvious since the
factor group $\brG_n / R\simeq \O_n /[\O_n ,\O_n]\simeq \GL_n (\Z
) / [\GL_n (\Z ), \GL_n (\Z )]\simeq \Z_2$ is abelian. Now, it is
obvious that $\brG_n /[\brG_n , \brG_n]\simeq \Z_2$.  $\Box $

$\noindent $

{\bf A characterization of the elements of the group $\rG_n$}. For
each automorphism $\s $ of the algebra $\mI_n$, the next lemma
gives explicitly the map $\v \in \Aut_K(P_n)$ such that $\s =
\s_\v$ (see Corollary \ref{d21Mar9}.(2)). Lemma  \ref{a18Apr9} is
used at the final stage of the proof of Theorem \ref{25Oct9}.

\begin{lemma}\label{a18Apr9}
For each automorphism $\s $ of the algebra $\mI_n$, there exists a
$K$-basis $\{ x'^{[\alpha]}\}_{\alpha \in \N^n}$ of the polynomial
algebra $P_n$ such that $\s (H_i) *x'^{[\alpha]} = (\alpha_i + 1)
x'^{[\alpha]}$ and $\s (\der_i) *x'^{[\alpha]} = x'^{[\alpha
-e_i]}$ for all $i=1, \ldots , n$ (where $x'^{[\beta]} :=0$ if
$\beta \in \Z^n\backslash \N^n$). Moreover,
\begin{enumerate}
\item  $\s =\s_\v$ where the map $\v \in \Aut_K(P_n)$:
$x^{[\alpha]} \mapsto x'^{[\alpha]}$ is the change-of-the-basis
map, \item $\s (\int_i) *x'^{[\alpha]} = x'^{[\alpha +e_i]}$ for
all $i=1, \ldots , n$, and \item  the basis $\{
x'^{[\alpha]}\}_{\alpha \in \N^n}$ is unique up to a simultaneous
multiplication of each element of the basis by the same nonzero
scalar.
\end{enumerate}
\end{lemma}

{\it Proof}.  Recall that the polynomial algebra
$P_n=\bigoplus_{\alpha \in \N^n}Kx^{[\alpha]}$ (where
$x^{[\alpha]}:=\prod_{i=1}^n \frac{x_i^{\alpha_i}}{i!}$) is the
direct sum of non-isomorphic, one-dimensional, simple $K[H_1,
\ldots , H_n]$-modules (see (\ref{kerHia})) such that
$\der_i*x^{[\alpha]} =x^{[\alpha - e_i]}$ for all $\alpha \in
\N^n$ and $i=1,\ldots , n$ (where $x^{[\beta]} :=0$ if $\beta \in
\Z^n\backslash \N^n$). Recall that $\s = \s_\v$ for some linear
map $\v \in \Aut_K(P_n)$, the linear map $\v : P_n\ra {}^\s P_n$
is an $\mI_n$-module isomorphism (Corollary \ref{d21Mar9}.(2)),
and the map $\v$ is unique up to a multiplication by a nonzero
scalar since $\End_{\mI_n} (P_n) \simeq K$ (Corollary
\ref{y16Apr9}.(1)).
 Let $x'^{[\alpha]} := \v (x^{[\alpha]} )$ for $\alpha\in \N^n$. Then the
 fact that the map $\v$ is an $\mI_n$-module homomorphism is equivalent to
 the fact that the following equations hold:
 \begin{eqnarray*}
 \s (H_i)*x'^{[\alpha]} &=& \v H_i\v^{-1}\v *x^{[\alpha]} = (\alpha_i+1) \v *x^{[\alpha]} =(\alpha_i+1) x'^{[\alpha]} ,  \\
  \s (\der_i)*x'^{[\alpha]} &=& \v \der_i\v^{-1}\v *x^{[\alpha]} =  \v *x^{[\alpha -e_i]} = x'^{[\alpha_i-e_i]} , \;\;\;\;\; (x'^{[\beta]}:=0, \; \beta \in \Z^n \backslash \N^n) \\
   \s (\int_i)*x'^{[\alpha]} &=& \v \int_i\v^{-1}\v *x^{[\alpha]} =  \v *x^{[\alpha +e_i]} =
   x'^{[\alpha_i+e_i]}.
\end{eqnarray*}
Note that the last equality follows from the previous two: by the
first equality, the polynomial algebra $P_n=\bigoplus_{\alpha \in
\N^n}Kx'^{[\alpha]}$ is the direct sum of non-isomorphic,
one-dimensional, simple $K[\s (H_1),\ldots , \s (H_n)]$-modules.
Since $\s (H_i) \s (\int_i) *x'^{[\alpha]} = \s
(H_i\int_i)*x'^{[\alpha]} = \s (\int_i(H_i+1))*x'^{[\alpha]}= \s
(\int_i) (\s (H_i) +1)*x'^{[\alpha]} = (\alpha_i+2) \s
(\int_i)*x'^{[\alpha]}$ for all $i$, we have $\s (\int_i)
*x'^{[\alpha]} = \l_{i,\alpha} x'^{[\alpha +e_i]}$ for a scalar $
\l_{i,\alpha}$ which is necessarily equal to $1$ since
$$ x'^{[\alpha]} = \s (\der_i)\s (\int_i)*x'^{[\alpha]} = \l_{i,\alpha} \s (\der_i)
*x'^{[\alpha +e_i]}= \l_{i,\alpha} x'^{[\alpha]}.$$ Since the
isomorphism $\v$ is unique up to a multiplication by a nonzero
scalar, the basis $\{ x'^{[\alpha]} \}$ is unique up to a
simultaneous multiplication of each element of it by the same
nonzero scalar. The proof of the lemma is complete.  $\Box $

$\noindent $

{\bf The group $\rG_n$ is a subgroup of $\mG_n$}.  In
\cite{algintdif}, it is proved that the Jacobian algebra $\mA_n=
S^{-1}\mI_n$ is the two-sided localization of the algebra $\mI_n$
at the multiplicatively closed subset $S:= \{
\prod_{i=1}^n(H_i+\alpha_i)_*^{n_i} \, | \, (\alpha_i) \in \Z^n,
(n_i)\in \N^n \}$ of $\mI_n$ where
$(H_i+\alpha_i)_*:=\begin{cases}
H_i+\alpha_i & \text{if }\alpha_i\geq 0,\\
(H_i+\alpha_i)_1& \text{if }\alpha_i<0,\\
\end{cases} $
and $(H_i-j)_1:= H_i-j+e_{j-1, j-1}(i)$ for $j\geq 1$. The
elements of the set $S\subseteq \End_K(P_n)$ are invertible linear
maps in $P_n$, i.e. $S\subseteq \Aut_K(P_n)$, and so are regular
elements of the algebra $\mI_n$ since $\mI_n\subseteq
\End_K(P_n)$.

\begin{theorem}\label{B25Oct9}
\begin{enumerate}
\item $\rG_n=\{ \s \in \mG_n \, | \, \s (\mI_n) = \mI_n\}$ and
$\rG_n$ is a subgroup of $\mG_n$. \item Each automorphism of the
algebra $\mI_n$ has a unique extension to an automorphism of the
algebra $\mI_n$.
\end{enumerate}
\end{theorem}

{\it Proof}. 1. Statement 1 follows from statement 2: the set $\{
\s \in \mG_n \, | \, \s (\mI_n) = \mI_n\}$ is a subgroup of the
group  $\mG_n$ that is mapped isomorphically onto the group
$\rG_n$ via $\s \mapsto \s|_{\mI_n}$, by statement 2.

2. It suffices to prove that each automorphism $\s$ of the algebra
$\mI_n$ can be extended to an automorphism of the algebra $\mA_n$,
 since then its uniqueness  is obvious as $\mA_n = S^{-1} \mI_n$. By
Corollary \ref{y16Apr9}.(2) and Lemma \ref{a18Apr9}, there exists
an $\mI_n$-module isomorphism $\v : P_n\ra {}^\s P_n$ which is
unique up to $K^*$, and  $\s (a) = \v a\v^{-1}$ for all $a\in
\mI_n$. Using the $K$-basis $\{ x'^{[\alpha]}\}_{\alpha \in \N^n}$
 of Lemma \ref{a18Apr9} we see that all the elements $\{ \s (s) \, |
\, s\in S\}$ are invertible in $\End_K({}^\s P_n) = \End_K(P_n)$.
By the universal property of localization, the algebra
monomorphism $\mI_n\ra \End_K({}^\s P_n)$, $ a\mapsto (p\mapsto \s
(a)p)$, can be extended uniquely to the algebra homomorphism
$\mA_n\ra \End_K({}^\s P_n)$, $s^{-1}a\mapsto (p\mapsto \s
(s)^{-1} \s (a) p)$ where $s\in S$ and $a\in \mI_n$. It is obvious
that the extension is an algebra monomorphism since $\s (S)
\subseteq \Aut_K({}^\s P_n)$. Therefore, the $\mA_n$-module ${}^\s
P_n$ is simple and faithful (since the $\mI_n$-module ${}^\s P_n$
is simple and $\mI_n\subseteq \mA_n$). The $\mA_n$-module $P_n$ is
 the only (up to isomorphism) simple and faithful $\mA_n$-module
(Corollary 2.7.(10), \cite{Bav-Jacalg}), and $\End_{\mA_n}(P_n)
\simeq K$. Therefore, there exists a unique (up to $K^*$)
$\mA_n$-module isomorphism $\psi : P_n\ra {}^\s P_n$  such that
$\s (a) = \psi a\psi^{-1}$ for all elements $a\in \mA_n$. In
particular, the map $\psi$ is an $\mI_n$-module isomorphism.
Therefore, $K^*\psi = K^*\v$ since $\End_{\mI_n}(P_n)\simeq K$.
Without loss of generality  we may assume that $\psi = \v$.
Therefore, the automorphism $\s \in \rG_n$ can be uniquely
extended to an automorphism of the algebra $\mA_n$. $\Box $

$\noindent $


For each natural number $d\geq 1 $, there is the decomposition
$K[x_i]= (\bigoplus_{j=0}^{d-1}Kx_i^j)\bigoplus (\bigoplus_{k\geq
d} Kx_i^k)$. The idempotents of the algebra $\mI_n$,
$p(i,d):=\sum_{j=0}^{d-1}e_{jj}(i)$ and $q(i,d):=1-p(i,d)$, are
the projections onto the first and the second summand
correspondingly. For a subset $I$ of the set $\{ 1, \ldots , n\}$,
 $CI$ denotes its complement.  Since $P_n=\bigotimes_{i=1}^n K[x_i]$, the identity map
$1=\id_{P_n}$ on the vector space  $P_n$ is the sum
\begin{equation}\label{1pIq}
1=\bigotimes_{i=1}^n (p(i,d)+q(i,d))=\sum_{I\subseteq \{ 1, \ldots
, n\}}p(I,d)q(CI,d)
\end{equation}
of orthogonal idempotents where $p(I,d):= \prod_{i\in I}p(i,d)$,
$q(CI,d):= \prod_{i\in CI}q(i,d)$, $p(\emptyset , d):=1$, and
 $q(\emptyset , d):=1$. Each idempotent $p(I,d)q(CI,d)\in
  \mI_n \subset \End_K(P_n)$ is the projection onto
 the summand $P_n(I, d)$ in the following  decomposition of the vector space $P_n$,
\begin{equation}\label{PnSI}
P_n=\bigoplus_{I\subseteq \{ 1, \ldots , n\} }P_n(I,d), \;\;
P_n(I,d) :=\bigoplus \{ Kx^{[\alpha]} \, | \, \alpha_i<d, \; {\rm
if}\; i\in I; \alpha_j\geq d, \; {\rm if}\; j\in CI\}.
\end{equation}
In particular, the idempotent $q(\{ 1, \ldots , n\}, d)$ is the
projection onto the subspace $P_n(\{ 1, \ldots , n\}, d)=\bigoplus
\{ Kx^{[\alpha]}\, | \, {\rm all}\; \alpha_i\geq d\}$.

$\noindent $

{\bf The group $\rG_n$ and  a formula for the map $\v $ such that
$\s = \s_\v$}. By Theorem \ref{21Mar9}, each element $\s = \s_\v
\in \rG_n$ (Corollary \ref{d21Mar9}.(2)) is uniquely determined by
the elements $\s (\der_1), \ldots , \s (\der_n)$.
 In the proof of Theorem \ref{25Oct9}, an explicit formula for the
 map $\v$ is given, (\ref{vf2Nov9}), via the elements $\s (\der_1),  \ldots , \s (\der_n)$.

By the very definition, the group $\ker (\xi )$ (see
(\ref{Ixiaut})) contains precisely all the automorphisms $\s \in
\rG_n$ such that 
\begin{equation}\label{sIdH}
 \s (\int_i) \equiv \int_i\mod \ga_n ,\;\;  \s (\der_i) \equiv
\der_i\mod \ga_n , \;\; \s (H_i) \equiv H_i\mod \ga_n ,\;\; i=1,
\ldots , n.
\end{equation}

\begin{theorem}\label{25Oct9}
\begin{enumerate}
\item $\rG_n = S_n \ltimes \mT^n  \ltimes \Inn (\mI_n )$. \item
$\rG_n = S_n \ltimes \mT^n \ltimes \ker (\xi ) $. \item $\im (\xi
) = S_n \ltimes \mT^n$. \item $\ker (\xi ) = \Inn (\mI_n)$.
\end{enumerate}
\end{theorem}

{\it Proof}. 1. Statement 1 follows from statements 2 and 4.

2. Statement 2 follows from statement 3: suppose that $\im (\xi )
= S_n \ltimes \mT^n$, then the  homomorphism $\xi$ maps
isomorphically the subgroup  $S_n\ltimes \mT^n $ of $\rG_n$ onto
its image, and so statement 2 follows from  the short exact
sequence of groups: $1\ra \ker (\xi ) \ra \rG_n \ra \im (\xi ) \ra
1$.

3. Let $\s \in \rG_n$. We have to show that there exists an
element $\s'\in S_n \ltimes \mT^n$ such that the automorphism
$\s'\s$ satisfies the conditions  (\ref{sIdH}). By Theorem
\ref{B25Oct9}, $\s \in \mG_n$. By (Corollary 7.5, \cite{jacaut}),
$\s (H_i) \equiv H_{\tau (i)}\mod \ga_n^e$ for all $i=1, \ldots ,
n$ and for some element $\tau \in S_n$ where $\ga_n^e$ is the only
maximal ideal of the algebra $\mA_n$ (Theorem \ref{7Oct9}). Then
$\tau^{-1}\s (H_i) \equiv H_i\mod \ga_n^e$ for all $i=1, \ldots ,
n$, and so $\tau^{-1}\s (H_i) - H_i\in \mI_n\cap \ga_n^e=
\ga_n^{er} = \ga_n$ (Theorem \ref{7Oct9}). For the automorphism
$\xi (\tau^{-1} \s ) \in \brG_n$, we have the action (\ref{atsp}),
$$ \xi (\tau^{-1} \s ) : z_i\mapsto \l_i z_i, \;\; H_i\mapsto H_i,
\;\; i=1, \ldots , n,$$ for some element $(\l_i) \in K^{*n}$. Then
$t_\l^{-1} \tau^{-1} \s \in \ker (\xi )$ where $\l = (\l_i)$.
Therefore, $\im (\xi ) = S_n\ltimes \mT^n$.

4. By Theorem \ref{A25Oct9}.(2), $\ker (\xi ) \supseteq \Inn
(\mI_n)$. Let $\s \in \ker (\xi )$. It remains to show that $\s
\in \Inn (\mI_n)$.  Fix a natural number $d$ such that
\begin{equation}\label{1sHdint}
\s (H_i)-H_i, \; \s (\der_i) - \der_i,\;  \s(\int_i) - \int_i\in
\sum_{k=1}^n \mI_{n-1, k}\t (\sum_{s,t=0}^{d-1} Ke_{st}(k))
\end{equation}
for all $i=1, \ldots , n$ where $\mI_{n-1,k}:= \bigotimes_{j\neq
k}\mI_1(j)$. By Lemma \ref{a18Apr9}, $\s = \s_\v:a\mapsto \v
a\v^{-1}$ where $\v \in \Aut_K(P_n):x^{[\alpha]}\mapsto
x'^{[\alpha]}$ is the change-of-the-basis map (see Lemma
\ref{a18Apr9}). By the choice of the number $d$ above, for each
element $x^{[\alpha]}\in P(\emptyset , d)$, $\s (H_i)
*x^{[\alpha]}= H_i*x^{[\alpha]}=(\alpha_i+1)x^{[\alpha]}$ and $\s
(\der_i) *x^{[\alpha]}=\der_i*x^{[\alpha]}  =x^{[\alpha-e_i]}$ for
all $i=1, \ldots , n$. By multiplying the map $\v$ by a nonzero
scalar, by Lemma \ref{a18Apr9}, we may assume that
\begin{equation}\label{xpa=xa}
x'^{[\alpha]}=x^{[\alpha]}\;\; {\rm for \; all}\;\; \alpha=
(\alpha_i) \in \N^n\;\; {\rm such \; that}\;\; \alpha_1\geq d,
\ldots , \alpha_n\geq d.
\end{equation}
It suffices to show that $\v \in \mI_n$ (since then $\v^{-1}\in
\mI_n$ as $\s^{-1} =\s_{\v^{-1}}$). This is obvious since

\begin{equation}\label{vf2Nov9}
\v = q(\{ 1, \ldots , n\}, d)+\sum_{\emptyset \neq I\subseteq  \{
1, \ldots , n\}}(\sum_{\alpha \in C_d(I)}\prod_{j\in I}\s
(\der_j^{d-\alpha_j})\cdot \prod_{i\in I}\int_i^{d-\alpha_i}\cdot
e_{\alpha \alpha}(I))p(I,d)q(CI, d)
\end{equation}
where $C_d(I):=\{ (\alpha_i)_{i\in I}\in \N^I\, | \, {\rm all }\;
\alpha_i<d\}$,  $e_{\alpha \alpha}(I):=\prod_{i\in
I}e_{\alpha_i\alpha_i}(i)$, and $d$ is as in  (\ref{1sHdint}). To
prove that this formula holds for the map  $\v$ we have to show
that $\v *x^{[\alpha]} = x'^{[\alpha]}$ for all $\alpha \in \N^n$.
For each $\alpha$, let $I:= \{ i \, | \, \alpha_i<d\}$. Then
$x^{[\alpha]} = \prod_{i\in I}x_i^{[\alpha_i]}\cdot \prod_{k\in
CI}x_k^{[\alpha_k]}$. If $I\neq \emptyset$ then
\begin{eqnarray*}
 \v *x^{[\alpha]} &=& \prod_{j\in I}\s (\der_j^{d-\alpha_j})\cdot
 \prod_{i\in I}\int_i^{d-\alpha_i}*\prod_{i\in
I}x_i^{[\alpha_i]}\cdot \prod_{k\in CI}x_k^{[\alpha_k]}=
\prod_{j\in I}\s (\der_j^{d-\alpha_j})*\prod_{i\in
I}x_i^{[d]}\cdot \prod_{k\in CI}x_k^{[\alpha_k]}  \\
&=& \prod_{j\in I}\s (\der_j^{d-\alpha_j})*x'^{[\sum_{i\in
I}de_i+\sum_{k\in CI}\alpha_ke_k]}=x'^{[\sum_{i\in
I}\alpha_ie_i+\sum_{k\in CI}\alpha_ke_k]}=x'^{[\alpha ]} .
\end{eqnarray*}
If $I=\emptyset$ then $\v *x^{[\alpha]} = q(\{ 1, \ldots , n\} ,
d) *x^{[\alpha]} = x^{[\alpha]}= x'^{[\alpha]}$.  The proof of the
theorem is complete. $\Box $

\begin{corollary}\label{a9Oct9}
Let $\s \in \Inn (\mI_n)$. Then there is a unique element $\v \in
(1+\ga_n)^*$ such that $\s (a) = \v a\v^{-1}$ for all elements
$a\in \mI_n$, and the element $\v $ is given by the formula
(\ref{vf2Nov9}).
\end{corollary}

{\it Proof}. The element $\v\in (1+\ga_n)^*$ such that $\s (a) =
\v a\v^{-1}$ for all $a\in \mI_n$ is unique by Theorem
\ref{A25Oct9}.(2). By the very definition, the element $\v \in
\mI_n^*$ from (\ref{vf2Nov9}) satisfies $\v \equiv 1\mod \ga_n$
and $\s (a) = \v a\v^{-1}$ for all $a\in \mI_n$.  Therefore, both
$\v$'s coincide. $\Box $

\begin{corollary}\label{b9Oct9}
$\Out (\mI_n) \simeq S_n \ltimes \mT^n$.
\end{corollary}

{\it Proof}. $\Out (\mI_n) =\rG_n / \Inn (\mI_n)= S_n \ltimes
\mT^n \ltimes \Inn (\mI_n)/\Inn (\mI_n)\simeq  S_n \ltimes \mT^n$.
$\Box$

$\noindent $

Recall that  $\CH_1:= \{ \gp_1, \ldots , \gp_n\}$ is the set of
height one prime ideals of the algebra $\mI_n$. The next corollary
describes its stabilizer $\St_{\rG_n}(\CH_1) := \{ \s \in \rG_n \,
| \, \s (\gp_1) = \gp_1, \ldots , \s (\gp_n) = \gp_n\}$.

\begin{corollary}\label{c9Oct9}
$\St_{\rG_n}(\CH_1) = \mT^n\ltimes \Inn (\mI_n)$.
\end{corollary}

{\it Proof}. It is obvious that $\St_{\rG_n}(\CH_1)\supseteq R :=
\mT^n\ltimes \Inn (\mI_n)$ and $S_n\cap \St_{\rG_n}(\CH_1) = \{
e\}$. Now,
$$ \St_{\rG_n}(\CH_1)= \rG_n \cap \St_{\rG_n}(\CH_1)= (S_n \ltimes
R)\cap \St_{\rG_n}(\CH_1)= (S_n \cap \St_{\rG_n}(\CH_1))\ltimes R
= R. \;\;\; \Box $$

The algebra $\mI_n = \bigoplus_{\alpha \in \Z^n}\mI_{n, \alpha}$
is a $\Z^n$-graded subalgebra of the Jacobian algebra $\mA_n =
\bigoplus_{\alpha \in \Z^n} \mA_{n,\alpha}$, see \cite{algintdif}.
 The group $\Aut_{\Z^n - {\rm gr}}(\mA_n) := \{ \s \in \mG_n \, |
\, \s (\mA_{n, \alpha } ) = \mA_{n, \alpha}$ for all $\alpha \in
\Z^n\}$ of $\Z^n$-grading preserving automorphisms  of the algebra
$\mA_n$ is equal to $\St_{\mG_n} (H_1, \ldots , H_n) :=\{ \s \in
\mG_n \, | \, \s (H_1) = H_1, \ldots , \s (H_n) = H_n \}$ and
$\St_{\mG_n}(H_1, \ldots , H_n) = \mT^n \times \mU_n$ (Corollary
7.10, \cite{jacaut}) where the subgroup $\mU_n$ of $\mG_n$ is
defined as follows. Recall that $\mA_n = S^{-1} \mI_n$. Let
$\CH_n$ be the subgroup of $\mA_n^*$ generated by the commutative
monoid $S\subseteq \mA_n^*$. Then the group $\CH_n=\prod_{i=1}^n
\CH_1(i)$ is the direct product of its subgroups $\CH_1(i) := \{
\prod_{j\geq 0}(H_i+j)^{n_j}\cdot \prod_{j\geq 1}
(H_i-j)^{n_{-j}}_1\, | \, (n_k)_{k\in \Z}\in \Z^{(\Z )}\} \simeq
\Z^{(\Z )}$, and so $\CH_n\simeq (\Z^n )^{(\Z )}$. Each element
$u=u_1\cdots u_n \in \CH_n$ where $u_i\in \CH_1(i)$ determines the
automorphism $\mu_u$ of the algebra $\mA_n$ (see \cite{jacaut} for
details),
$$ \mu_u:x_i\mapsto x_iu_i, \;\; y_i\mapsto u_i^{-1}y_i, \;\;
H_i^{\pm 1}\mapsto H_i^{\pm 1}, \;\; i=1, \ldots , n.$$ Then
$\mU_n := \{ \mu_u\, | \, u\in \CH_n\}\simeq (\Z^n )^{(\Z )}$. Let
 $\St_{\rG_n} (H_1, \ldots , H_n):=\{ \s \in \rG_n\, | \, \s (H_1)
 = H_1, \ldots , \s (H_n) = H_n\}$.

\begin{corollary}\label{d9Oct9}
\begin{enumerate}
\item $\St_{\rG_n} (H_1, \ldots , H_n)=\mT^n$. \item Let $\s ,
\tau \in \rG_n$. Then  $\s (H_1) =\tau (H_1) ,\ldots , \s (H_n) =
\tau (H_n)$ iff $\s = \tau t_\l$ for some element $t_\l \in
\mT^n$. \item $\Aut_{\Z^n - {\rm gr}}(\mI_n)= \St_{\rG_n} (H_1,
\ldots , H_n)\subset \Aut_{\Z^n - {\rm gr}}(\mA_n)$.
\end{enumerate}
\end{corollary}

{\it Proof}. 1. Since $\rG_n \subseteq \mG_n$ (Theorem
\ref{B25Oct9}),
$$\St_{\rG_n} (H_1, \ldots , H_n)=\rG_n\cap \St_{\mG_n} (H_1, \ldots ,
H_n)=\rG_n\cap \mT^n \times \mU_n= \mT^n \times (\rG_n \cap \mU_n)
= \mT^n$$ since $\rG_n\cap \mU_n = \{ e\}$, by the very definition
of the group $\mU_n$.

2. Statement 2 follows from statement 1.

3. Since $\mA_n = S^{-1} \mI_n$, $S\subseteq \mI_{n,0}$ and $\rG_n
\subseteq \mG_n$, we have the inclusion $\Aut_{\Z^n - {\rm
gr}}(\mI_n)\subseteq \Aut_{\Z^n - {\rm gr}}(\mA_n)$. Now,
$\Aut_{\Z^n - {\rm gr}}(\mI_n)=\rG_n\cap \Aut_{\Z^n - {\rm
gr}}(\mA_n)= \rG_n \cap \mT^n \times \mU_n = \mT^n =\St_{\rG_n}
(H_1, \ldots , H_n)$, by statement 1. By statement 1, the
inclusion $\Aut_{\Z^n - {\rm gr}}(\mI_n)\subseteq \Aut_{\Z^n -
{\rm gr}}(\mA_n)$ is a strict inclusion  since $\Aut_{\Z^n - {\rm
gr}}(\mA_n)=\mT^n \times \mU_n$ (Corollary 7.10, \cite{jacaut}).
 $\Box $

$\noindent $

{\bf The canonical form of $\s \in \rG_n$}. By Theorem
\ref{25Oct9}, each automorphism $\s$ of the algebra $\mI_n$ is a
{\em unique} product $st_\l \o_\v$ where $s\in S_n$, $t_\l \in
\mT^n$,  and $\o_\v$ is an inner automorphism of the algebra
$\mI_n$ with $\v \in (1+\ga_n)^*$, and the element $\v$ is unique
(Corollary \ref{a9Oct9}).

$\noindent $

{\it Definition}. The unique product $\s =st_\l \o_\v$  is called
the {\em canonical form} of the automorphism $\s $ of the algebra
$\mI_n$.

\begin{corollary}\label{a2Nov9}
Let $\s \in \rG_n$ and $ \s =st_\l\o_\v$ be its canonical form.
Then the automorphisms $s$, $t_\l$  and $\o_\v$ can be effectively
(in finitely many steps) found from the action of the automorphism
$\s$ on the elements $\{ H_i, \der_i, \int_i\, | \, i=1, \ldots ,
n\}$:
$$ \s (H_i) \equiv H_{s(i)}\mod \ga_n, \;\; \s (\der_i) \equiv
\l_i^{-1} \der_{s_(i)}\mod \ga_n,\;\; \s (\int_i) \equiv
\l_i\int_{s(i)}\mod \ga_n,
$$ and the elements $\v$  and $\v^{-1}$ are  given by the formulae
 (\ref{vf2Nov9}) and (\ref{1vf2Nov9}) respectively
for the automorphism $(st_\l)^{-1}\s \in \Inn (\mI_n)=\ker (\xi
)$.
\end{corollary}
The next corollary is a criterion for an automorphism of the
algebra $\mI_n$ to be an inner automorphism.
\begin{corollary}\label{b2Nov9}
Let $\s \in \rG_n$. The following statements are equivalent.
\begin{enumerate}
\item $\s \in \Inn (\mI_n)$. \item $\s (\der_i) \equiv \der_i \mod
\ga_n$ for $i=1, \ldots , n$. \item $\s (\int_i) \equiv \int_i
\mod \ga_n$ for $i=1, \ldots , n$.
\end{enumerate}
\end{corollary}

{\it Proof}.  The result follows from Theorem \ref{25Oct9}.(4) and
Corollary \ref{a2Nov9}. $\Box $

\begin{corollary}\label{c2Nov9}
Let $\s \in \rG_n$. Then $\s \in \mT^n \ltimes \Inn (\mI_n)$ iff
$\s (H_i) \equiv H_i\mod \ga_n$ for $i=1, \ldots , n$.
\end{corollary}

{\it Proof}. This follows from  Theorem \ref{25Oct9}.(4) and
Corollary \ref{a2Nov9}. $\Box $

$\noindent $

{\bf A formula for the inverse $\v^{-1}$ where $\s= \s_\v \in \Inn
(\mI_n)$ via $\s (\der_i)$ and $\s (\int_j )$}. By Corollary
\ref{a9Oct9}, for each inner automorphism $\s \in \Inn (\mI_n)$
there exists a unique element $\v\in (1+\ga_n)^*$ such that $\s =
\s_\v : a\mapsto \v a\v^{-1}$ for all $a\in \mI_n$. The next
theorem presents a formula for the inverse $\v^{-1}$ via the
elements $\{ \s (\der_i), \s (\int_i)\, | \, i=1, \ldots , n\}$.

\begin{theorem}\label{12Nov9}
Let $\s = \s_\v\in \Inn (\mI_n)$ where $\v \in (1+\ga_n)^*$
($\s_\v (a) = \v a\v^{-1}$ for all $a\in \mI_n$, see Corollary
\ref{a9Oct9}). Then $\s_\v^{-1} = \s_{\v^{-1}}$ and
\begin{equation}\label{1vf2Nov9}
\v^{-1} = q'(\{ 1, \ldots , n\}, d)+\sum_{\emptyset \neq
I\subseteq \{ 1, \ldots , n\}}(\sum_{\alpha \in C_d(I)}\prod_{j\in
I}\der_j^{d-\alpha_j}\cdot \prod_{i\in I}\s
(\int_i)^{d-\alpha_i}\cdot e_{\alpha \alpha}'(I))p'(I,d)q'(CI, d)
\end{equation}
where $d$ is as in (\ref{1sHdint}) for $\s = \s_\v$,  $e_{\alpha
\alpha}'(I):=\prod_{i\in I}e_{\alpha_i\alpha_i}'(i)$ and
$e_{jj}'(i):=\s (e_{jj}(i))=\s (\int_i)^j\s (\der_i)^j-\s
(\int_i)^{j+1}\s (\der_i)^{j+1}$; $p'(I,d):=\prod_{i\in I}p'(i,d)$
and $p'(i,d):=  \s (p(i,d))=\sum_{j=0}^{d-1}e_{jj}'(i)$;  $q'(CI,
d):= \s (q(CI, d))=\prod_{i\in CI}(1-p'(d,i))$.
\end{theorem}

{\it Proof}. We keep the notation of the proof of statement 4 of
Theorem \ref{25Oct9}. In particular, $\v : P_n = \bigoplus_{\alpha
\in \N^n} Kx^{[\alpha]}\ra P_n=\bigoplus_{\alpha \in \N^n}
Kx'^{[\alpha]}$, $x^{[\alpha]}\mapsto x'^{[\alpha]}$. For each
$\alpha \in \N^n$, the projection onto the summand
$Kx'^{[\alpha]}$ of the polynomial algebra $P_n$ is equal to $\v
e_{\alpha\alpha}\v^{-1}= \s (e_{\alpha\alpha})$. For each subset
$I$ of the set $\{ 1, \ldots , n\}$, let
$$P_n'(I,d):=\v (P_n(I,d))=\bigoplus \{ Kx'^{[\alpha]} \, | \,
 \alpha_i<d\; {\rm if}\; i\in I; \alpha_j\geq d\; {\rm if}\; j\in
 CI\},$$ see (\ref{PnSI}).  Since
 $\v : P_n=\bigoplus_{I\subseteq \{ 1, \ldots , n\}} P_n(I,d)
\simeq P_n=\bigoplus_{I\subseteq \{ 1, \ldots , n\}} P_n'(I,d)$,
the projections onto the summand $P'_n(I,d)$ of $P_n$ is equal to
$$ \v p(I,d)q(I,d) \v^{-1} = \s (p(I,d)q(I,d)) = \s (p(I,d)) \s
(q(I,d))= p'(I,d)q'(I,d)$$ where $p'(I,d)= \s (p(I,d))$ and
$q'(I,d)=\s (q(I,d))$. Then the inverse map $\v^{-1}$ of the map
$\v$ in (\ref{vf2Nov9}) is given by (\ref{1vf2Nov9}).  To prove
this let $\psi$ be the RHS of (\ref{1vf2Nov9}). We have to show
that $\psi : x'^{[\alpha]}\mapsto x^{[\alpha]}$  for all $ \alpha
\in \N^n$. Fix $\alpha$, and let $I:= \{ i \, | \, \alpha_i<d\}$.
Then $x'^{[\alpha]} = \prod_{i\in I}x_i'^{[\alpha_i]}\cdot
\prod_{k\in CI}x_k'^{[\alpha_k]}$. If $I\neq \emptyset$ then,  by
Lemma \ref{a18Apr9},
\begin{eqnarray*}
 \psi *x'^{[\alpha]} &=& \prod_{j\in I}\der_j^{d-\alpha_j}\cdot \prod_{i\in I}\s (\int_i)^{d-\alpha_i}*\prod_{i\in
I}x_i'^{[\alpha_i]}\cdot \prod_{k\in CI}x_k'^{[\alpha_k]}=
\prod_{j\in I}\der_j^{d-\alpha_j}*\prod_{i\in
I}x_i'^{[d]}\cdot \prod_{k\in CI}x_k'^{[\alpha_k]}  \\
&=& \prod_{j\in I}\der_j^{d-\alpha_j}*\prod_{i\in I}x_i^{[d]}\cdot
\prod_{k\in CI}x_k^{[\alpha_k]}\;\;\;\; ({\rm by}\;\; (\ref{xpa=xa}))\\
&=&x^{[\alpha]}.
\end{eqnarray*}
If $I\neq \emptyset$ then $\psi *x'^{[\alpha]} = q'(\{ 1, \ldots ,
n\} , d)*x'^{[\alpha]} = x'^{[\alpha]} = x^{[\alpha]}$, by
(\ref{xpa=xa}). This finishes the proof of the theorem. $\Box $

$\noindent $

{\bf An inversion formula for $\s \in \rG_n$}. The next theorem
gives an inversion formula for $\s$ via the elements $\{ \s
(\der_i) , \s(\int_i)\, | \, i=1, \ldots , n\}$.

\begin{theorem}\label{13Nov9}
Let $\s \in \rG_n$ and $\s = st_\l \o_\v$ be its canonical form
where $s\in S_n$,  $t_\l\in \mT^n$ and $\o_\v \in \Inn (\mI_n)$
for a unique element $\v \in (1+\ga_n)^*$. Then 
\begin{equation}\label{rGinv}
\s^{-1} = s^{-1} t_{s(\l^{-1})}\o_{st_\l (\v^{-1})}
\end{equation}
is the canonical form of the automorphism $\s^{-1}$ where the
elements  $\v^{-1}$ and $\v$  are given by
 the formulae (\ref{1vf2Nov9}) and (\ref{vf2Nov9}) respectively  for the automorphism
 $(st_\l )^{-1}\s \in \Inn (\mI_n)$.
\end{theorem}

{\it Proof}. $\s^{-1} = s^{-1} \cdot st^{-1}_\l s^{-1} \cdot st_\l
\o_\v^{-1} (st_\l)^{-1}=s^{-1} t_{s(\l^{-1})}\o_{st_\l
(\v^{-1})}$. $\Box $

\begin{theorem}\label{14Nov9}
The centre of the group $\rG_n$ is $\{ e\}$.
\end{theorem}

{\it Proof}. Let $\s $ be an element of the centre of the group
$\rG_n$. For all elements $\alpha , \beta \in \N^n$,
$1+e_{\alpha\beta}\in (1+\ga_n)^*$, and so
$\o_{1+e_{\alpha\beta}}\in \Inn (\mI_n)$. Then
$\o_{1+e_{\alpha\beta}}=\s \o_{1+e_{\alpha\beta}}\s^{-1}= \o_{1+\s
(e_{\alpha\beta})}$, and so $1+e_{\alpha\beta}=1+\s
(e_{\alpha\beta})$ (Theorem \ref{A25Oct9}.(2)), i.e.
$e_{\alpha\beta} = \s (e_{\alpha\beta})$. By Corollary
\ref{yz16Apr9}, $\s =e$.  $\Box $

$\noindent $

Let $H$ be a subgroup of a group $G$. The {\em centralizer}
$\Cen_G(H):=\{ g\in G\, | \, gh=hg\;\; {\rm for \; all}\; h\in
H\}$ of $H$ in $G$ is a subgroup of $G$. In the proof of Theorem
\ref{14Nov9}, we have used only inner derivations of the algebra
$\mA_n$. So, in fact, we have proved there  the next corollary.

\begin{corollary}\label{a14Nov9}
$\Cen_{\rG_n}(\Inn (\mI_n))=\{e\}$.
\end{corollary}

For an algebra $A$ and a subgroup $G$ of its group of algebra
automorphisms, the set $A^G:=\{ a\in A\, | \, \s (a) = a$ for all
$\s \in G\}$ is called the {\em fixed algebra} or the {\em algebra
of invariants}  for the group $G$.

\begin{theorem}\label{A14Nov9}
$\mI_n^{\rG_n}= \mI_n^{\Inn (\mI_n)}=K$.
\end{theorem}

{\it Proof}. Since $K\subseteq \mI_n^{\rG_n}\subseteq  \mI_n^{\Inn
(\mI_n)}$, it suffices to show that $\mI_n^{\Inn (\mI_n)}=K$. For
all $\alpha , \beta \in \N^n$, $1+e_{\alpha\beta}\in (1+\ga_n)^*$.
Then $\o_{1+e_{\alpha\beta}}\in \Inn (\mI_n)$. If $a\in\mI_n^{\Inn
(\mI_n)}$ then $a=\o_{1+e_{\alpha\beta}}(a)$, and so
$ae_{\alpha\beta} = e_{\alpha\beta} a$. By Lemma
\ref{c14Nov9}.(1), $a\in K$. Therefore, $\mI_n^{\Inn (\mI_n)}=K$.
 $\Box $

\begin{lemma}\label{c14Nov9}
\begin{enumerate}
\item $\Cen_{\mI_n}(\{ e_{\alpha\beta}\, | \, \alpha, \beta \in
\N^n\}) = K$. \item $\Cen_{\mI_n}(\ga ) = K$ for all nonzero
ideals $\ga$ of the algebra $\mI_n$.
\end{enumerate}
\end{lemma}

{\it Proof}. 1. Recall that $\mI_n\simeq \mI_n^{op}$ and $F_n
=\bigoplus_{\alpha, \beta \in \N^n}Ke_{\alpha  \beta}$ is the
least nonzero ideal of the algebra $\mI_n$. Hence, $F_n$ is a
simple $\mI_n$-bimodule. Since
$$ {}_{\mI_n}{F_n}_{\mI_n}\simeq {}_{\mI_n\t\mI_n^{op}}F_n\simeq  {}_{\mI_n\t\mI_n}F_n\simeq
 {}_{\mI_{2n}}F_n\simeq {}_{\mI_{2n}}P_{2n},$$
we have $K\simeq \End_{\mI_{2n}}(P_{2n})\simeq \Cen_{\mI_n}(F_n)
=\Cen_{\mI_n}(\{ e_{\alpha\beta}\, | \, \alpha, \beta \in
\N^n\})$.

2. Since $F_n\subseteq \ga$, we have $K\subseteq \Cen_{\mI_n}(\ga)
\subseteq \Cen_{\mI_n}(F_n)=K$, by statement 1. Therefore,
$\Cen_{\mI_n}(\ga ) = K$.  $\Box $

\begin{theorem}\label{B14Nov9}
No proper prime factor algebra of $\mI_n$ can be embedded into
$\mI_n$ (that is, for each nonzero prime ideal $\gp$ of the
algebra $\mI_n$, there is no algebra monomorphism from $\mI_n/\gp$
into $\mI_n$).
\end{theorem}

{\it Proof}. By Corollary \ref{b10Oct9}.(7), $\gp=
\gp_{i_1}+\cdots + \gp_{i_s}$. Without loss of generality, we may
assume that $\gp= \gp_1+\cdots + \gp_s$. Suppose that there is a
monomorphism $f:\mI_n/ \gp \ra \mI_n$, we seek a contradiction.
For each element $a\in \mI_n$, let $\oa := a+\gp$. Notice that,
for $i=1, \ldots ,  s$, $\overline{\der}_i\overline{\int}_i = 1$
and $\overline{\int}_i\overline{\der}_i =
\overline{1-e_{00}(i)}=1$ since $e_{00}(i)\in \gp_i \subseteq
\gp$. The elements $\{ \overline{\der}_i, \overline{\int}_i \, |
\, i=1, \ldots , s\}$ are units of the algebra $\mI_n / \gp$,
hence their images under the map $f$ are units of the algebra
$\mI_n$, i.e. $f(\int_i) , f(\der_i) \in \mI_n^* = K^* \times
(1+\ga_n)^*$ (Theorem \ref{A25Oct9}.(1)). We see that the image of
the {\em simple noncommutative} algebra $B_s:= \mI_s/ \ga_s$ under
the compositions of homomorphisms $B_s:= \mI_s/\ga_s \ra \mI_n /
\gp \stackrel{f}{\ra} \mI_n \ra B_n = \mI_n / \ga_n $ is the
subalgebra of $B_n$ generated by the images of the {\em
commutative} elements $H_1, \ldots , H_n$, a contradiction.  $\Box
$


$\noindent $

The next lemma shows that in the algebra $\mI_n$ there are
non-invertible elements that are invertible as elements of the
algebra $\End_K(P_n)$.

\begin{lemma}\label{b14Nov9}
$\mI_n^* \subsetneqq \mI_n \cap \Aut_K(P_n)$.
\end{lemma}

{\it Proof}. The element $1-\der_i\in \End_K(K[x_i])$ is an
invertible linear map since $\der_i$ is a locally nilpotent
derivation of the polynomial algebra $K[x_i]$, hence
$u:=\prod_{i=1}^n (1-\der_i) \in \mI_n\cap \Aut_K(P_n)$. But
$u\not\in \mI_n^*$ since the element $u+\ga_n$ is not a unit of
the factor algebra $B_n=\mI_n/ \ga_n$ as $B_n^* = \bigcup_{\alpha
\in \Z^n} K^*\der^\alpha$. Therefore, $\mI_n^* \subsetneqq \mI_n
\cap \Aut_K(P_n)$.  $\Box $

$\noindent $

{\bf The group $\CG_n := \Aut_{K-{\rm alg}}(\CI_n)$}.

{\it Definition}, \cite{shrekalg}. The 
{\em algebra} $\mathbb{S}_n$ {\em of one-sided inverses} of $P_n$
is an algebra generated over a field $K$ of arbitrary
characteristic by $2n$ elements $x_1, \ldots , x_n, y_n, \ldots ,
y_n$ that satisfy the defining relations:
$$ y_1x_1=\cdots = y_nx_n=1 , \;\; [x_i, y_j]=[x_i, x_j]= [y_i,y_j]=0
\;\; {\rm for\; all}\; i\neq j,$$ where $[a,b]:= ab-ba$ is  the
algebra  commutator of elements $a$ and $b$. Let $G_n:=
\Aut_{K-{\rm alg}}(\mS_n)$.

$\noindent $

By the very definition, the algebra $\mS_n\simeq \mS_1^{\t n}$ is
obtained from the polynomial algebra $P_n$ by adding commuting,
left (but not two-sided) inverses of its canonical generators. The
algebra $\mS_1$ is a well-known primitive algebra
\cite{Jacobson-StrRing}, p. 35, Example 2. Over the field
 $\mathbb{C}$ of complex numbers, the completion of the algebra
 $\mS_1$ is the {\em Toeplitz algebra} which is the
 $\mathbb{C}^*$-algebra generated by a unilateral shift on the
 Hilbert space $l^2(\N )$ (note that $y_1=x_1^*$). The Toeplitz
 algebra is the universal $\mathbb{C}^*$-algebra generated by a
 proper isometry.

The algebra $\CI_n := K\langle \der_1, \ldots , \der_n , \int_1,
\ldots , \int_n\rangle$ of scalar integro-differential operators
is isomorphic to the algebra $\mS_n$: 
\begin{equation}\label{SnIna}
 \mS_n \ra \CI_n, \;\; x_i\mapsto \int_i, \;\; y_i\mapsto \der_i,
\;\; i=1, \ldots , n.
\end{equation}

 Since $\CI_n =
\bigotimes_{i=1}^n\CI_1(i)$ where $\CI_1(i) := K\langle \der_i,
\int_i\rangle$ and $\mS_n =\bigotimes_{i=1}^n \mS_1(i)$ where
$\mS_1(i) := K\langle x_i, y_i\rangle$, it suffices to prove  the
statement for $n=1$. For $n=1$, the algebra epimorphism $\mS_1\ra
\mI_1$ is an isomorphism since any proper epimorphic image of the
algebra $\mS_1$ is commutative (see \cite{shrekalg}) but the
algebra $\CI_1$ is non-commutative. The algebra $\mS_n$ was
studied in detail in \cite{shrekalg},  its group
 of automorphism and explicit generators were  found in the
papers \cite{shrekaut}, \cite{K1aut} and  \cite{Snaut}.

\begin{theorem}\label{2Nov9}
\begin{enumerate}
\item $\CG_n = S_n \ltimes \mT^n \ltimes \Inn (\CI_n)$ and $\Inn
(\CI_n) =\{ \o_u\, | \, u\in (1+\ga_n')^*\} \simeq (1+\ga_n')^*$,
$\o_u\mapsto u$, where $\ga_n':= \sum_{i=1}^n \CI_n F(i)$ is the
only maximal ideal of the algebra $\CI_n$. \item $\CG_n =\{ \s \in
\rG_n \, | \, \s (\CI_n) = \CI_n \} = \{ \s \in \mG_n \, | \, \s
(\CI_n ) = \CI_n \}$ and $\CG_n$ is a subgroup of the groups
$\rG_n $ and $\mG_n$. \item Each automorphism of the algebra
$\CI_n$ has a unique extension to an automorphism of the algebra
$\mI_n$ and $\mA_n$. \item $\rG_n \supseteq \CG_n \supseteq
S_n\ltimes \mT^n\ltimes \underbrace{\GL_\infty (K)\ltimes\cdots
\ltimes \GL_\infty (K)}_{2^n-1 \;\; {\rm times}}$.
\end{enumerate}
\end{theorem}

{\it Proof}. 1. $\CI_n \simeq \mS_n$, $G_n = S_n \ltimes \mT^n
\ltimes \Inn (\mS_n)$ and $\Inn (\mS_n) =\{ \o_u\, | \, u\in
(1+\ga_n'')^*\} \simeq (1+\ga_n'')^*$, $\o_u\mapsto u$, where
$\ga_n''$ is the only maximal ideal of the algebra $\mS_n$,
\cite{shrekaut}.

2. Statement 2 follows from statement 1, Theorem \ref{25Oct9},
Theorem \ref{B25Oct9} and Corollary \ref{d21Mar9}.

3. Statement 3 follows from Theorem \ref{21Mar9}.

4. $G_n \supseteq S_n\ltimes \mT^n\ltimes \underbrace{\GL_\infty
(K)\ltimes\cdots \ltimes \GL_\infty (K)}_{2^n-1 \;\; {\rm
times}}$, \cite{shrekaut}. $\Box $



\section{Stabilizers  of the ideals of $\mI_n$ in $\rG_n$}\label{ISTABAN}

In this section, for each nonzero  ideal $\ga $ of the algebra
$\mI_n$ its stabilizer $\St_{G_n}(\ga ) :=\{ \s \in \rG_n \, | \,
\s (\ga )=\ga \}$ is found (Theorem \ref{15Apr9}) and it is  shown
that the stabilizer $\St_{\rG_n}(\ga )$ has finite index in the
group $\rG_n$ (Corollary \ref{d15Apr9}). When the ideal $\ga$ is
either prime or generic, this result can be refined even further
(Corollary \ref{b15Apr9}, Corollary \ref{c15Apr9}).
 In particular, when $n>1$ the stabilizer
of each height 1 prime  of $\mI_n$ is a maximal subgroup of
$\rG_n$ of index $n$ (Corollary \ref{b15Apr9}.(1)). It is shown
that the ideal $\ga_n$ is the only nonzero, prime,
$\rG_n$-invariant ideal of the algebra $\mI_n$ (Corollary
\ref{b15Apr9}.(3)).

 An ideal $\ga$  of
 $\mI_n$ is called  a {\em proper}
ideal if  $\ga\neq 0, \mI_n$. For an ideal $\ga$ of the algebra
$\mI_n$, $\Min (\ga )$ denotes  the set of all the minimal primes
over $\ga$. Two ideals $\ga$ and $\gb$ are called {\em
incomparable} if neither $\ga \subseteq \gb$ nor $\gb \subseteq
\ga$. The
 ideals of the algebra $\mI_n$ are classified  in
\cite{algintdif}. The next theorem shows that each ideal of the
algebra $\mI_n$ is completely determined by its minimal primes. We
use this theorem in the proof of Theorem \ref{15Apr9}.

\begin{theorem}\label{a15Apr9}
{\rm \cite{algintdif}} Let $\ga$ be a proper
 ideal of the algebra $\mI_n$. Then $\Min (\ga )$ is a
finite non-empty set, and the ideal $\ga$ is a unique product and
a unique intersection of incomparable prime  ideals of $\mI_n$
(uniqueness is  up to permutation). Moreover,
$$ a= \prod_{\gp \in \Min (\ga )}\gp =  \bigcap_{\gp \in \Min (\ga )}\gp .$$
\end{theorem}

Let $\Sub_n$ be the set of all the subsets of the set $\{ 1,
\ldots, n \}$. $\Sub_n$ is a partially ordered set with respect to
`$\subseteq $'. Let $\SSub_n$ be the set of all the  subsets of
$\Sub_n$. An element $\{ X_1, \ldots , X_s\}$ of $\SSub_n$ is
called an {\em antichain} if for all $i\ne j$ such that $1\leq
i,j\leq s$ neither $X_i\subseteq X_j$ nor $X_i\supseteq X_j$. An
empty set and one element set are antichains  by definition. Let
$\Inc_n$ be the subset of $\SSub_n$ that contains  all the
antichains of $\SSub_n$. The number $\gd_n:= |\Inc_n|$ is called
the {\em Dedekind} number. The symmetric group $S_n$ acts in the
obvious way on the sets $\SSub_n$ and $\Inc_n$ ($\s \cdot \{ X_1,
\ldots ,X_s\}= \{ \s (X_1), \ldots , \s (X_s)\}$).

\begin{theorem}\label{15Apr9}
Let $\ga $ be a proper  ideal of the algebra $\mI_n$. Then
$$ \St_{\rG_n}(\ga ) = \St_{S_n}(\Min (\ga ))\ltimes \mT^n\ltimes \Inn (\mI_n)
$$ where $\St_{S_n}(\Min (\ga )):=\{ \s \in S_n \, | \,
\s (\gq ) \in \Min (\ga )$ for all $\gq \in \Min (\ga )\}$.
Moreover, if $\Min (\ga ) = \{ \gq_1, \ldots , \gq_s\}$ and, for
each number $t=1, \ldots , s$, $\gq_t=\sum_{i\in I_t}\gp_i$ for
some subset $I_t$ of $\{ 1, \ldots , n\}$. Then the group
$\St_{S_n}(\Min (\ga ))$ is the stabilizer in the group $S_n$  of
the element $\{ I_1, \ldots , I_s\}$ of $\SSub_n$.
\end{theorem}

{\it Remark}. Note that the group $$\St_{\rG_n}(\Min (\ga ))=
\St_{S_n}(\{ I_1, \ldots , I_s\}):=\{ \s \in S_n \, | \, \{ \s
(I_1), \ldots , \s (I_s)\} = \{ I_1, \ldots , I_s\}\}$$ (and also
the group $\St_{\rG_n}(\ga )$) can be effectively computed in
finitely many steps.

$\noindent $

{\it Proof}. Recall that each nonzero prime ideal of the algebra
 $\mI_n$ is a unique sum of height one  prime ideals of the algebra
 $\mI_n$. By Theorem \ref{a15Apr9} and Corollary \ref{c9Oct9}, $\St_{\rG_n}(\ga )\supseteq \St_{\rG_n}(\CH_1
 )= \mT^n\ltimes
\Inn (\mI_n)$. Since $\rG_n=S_n\ltimes \mT^n\ltimes \Inn (\mI_n)$
(Theorem \ref{25Oct9}.(1)),
$$
\St_{\rG_n}(\ga ) =(\St_{\rG_n}(\ga )\cap S_n)\ltimes \mT^n\ltimes
\Inn (\mI_n)=\St_{S_n}(\ga )\ltimes \mT^n\ltimes\Inn (\mI_n).$$
 By
Theorem \ref{a15Apr9}, $\St_{S_n}(\ga )=\St_{S_n}(\Min (\ga ))=
\St_{S_n}(\{ I_1, \ldots , I_s\} )$, and the statement follows.
$\Box$

$\noindent $

The {\em index} of a subgroup $H$ in a group $G$ is denoted by
$[G:H]$.
\begin{corollary}\label{d15Apr9}
Let $\ga $ be a proper  ideal of $\mI_n$. Then $
[\rG_n:\St_{\rG_n}(\ga )] = |S_n:\St_{S_n}(\Min (\ga ))|<\infty$.
\end{corollary}

{\it Proof}. This follows from Theorem \ref{25Oct9}.(1) and
Theorem \ref{15Apr9} $\Box $

\begin{corollary}\label{b15Apr9}
\begin{enumerate}
\item $\St_{\rG_n}(\gp_i) \simeq S_{n-1}\ltimes \mT^n\ltimes \Inn
(\mI_n)$, for $i=1, \ldots , n$. Moreover, if $n>1$ then the
groups $\St_{\rG_n}(\gp_i)$ are maximal subgroups of $\rG_n$ with
$[\rG_n :\St_{\rG_n}(\gp_i)]=n$ (if $n=1$ then
$\St_{\rG_1}(\gp_1)=\rG_1$, see statement 3). \item Let $\gp$ be a
nonzero  prime ideal of the algebra $\mI_n$ and $h= \hht (\gp )$
be its height. Then $\St_{\rG_n}(\gp ) \simeq (S_h\times
S_{n-h})\ltimes  \mT^n\ltimes \Inn (\mI_n) $.\item The ideal
$\ga_n$ is the only nonzero, prime, $\rG_n$-invariant ideal of the
algebra $\mI_n$. \item Suppose that $n>1$. Let $\gp$ be a nonzero
prime ideal of the algebra $\mI_n$. Then its stabilizer
$\St_{\rG_n}(\gp )$ is a maximal subgroup of $\rG_n$ iff the ideal
$\gp $ is of height one.
\end{enumerate}
\end{corollary}

{\it Proof}. 1. Clearly, $\St_{\rG_n}(\gp_i) \cap S_n=\{ \tau \in
S_n\, | \, \tau (\gp_i) = \gp_i\} \simeq S_{n-1}$. By Theorem
\ref{15Apr9}, $\St_{\rG_n}(\gp_i)=S_{n-1}\ltimes \mT^n \ltimes
\Inn (\mI_n)$. When $n>1$, the group $\St_{\rG_n}(\gp_i)$ is a
maximal subgroup of $\rG_n$ since
$$S_{n-1}\simeq \St_{G_n}(\gp_i)/\mT^n \ltimes \Inn (\mI_n)
 \subseteq G_n/ \mT^n \ltimes \Inn (\mI_n)
\simeq S_n$$ and $S_{n-1}=\{ \s \in S_n \, | \, \s (i)=i \}$ is a
maximal subgroup of $S_n$. Clearly, $[\rG_n
:\St_{\rG_n}(\gp_i)]=[S_n:S_{n-1}]=n$.

2. By Corollary \ref{b10Oct9}.(9), $\gp=\gp_{i_1}+\cdots +
\gp_{i_h}$ for some distinct indices $i_1, \ldots , i_h\in \{ 1,
\ldots , n\}$. Let $I=\{ i_1, \ldots , i_h\}$ and $CI$ be its
complement in the set $\{ 1, \ldots , n\}$. Statement 2 follows
from Theorem \ref{15Apr9} and the fact that
$$\St_{\rG_n} (\gp ) \cap S_n = \{ \s \in S_n \, | \, \s (I) = I, \s
(CI) = CI\}\simeq S_h\times S_{n-h}.$$

3. Since $\ga_n = \gp_1+\cdots +\gp_n$, statement 3 follows from
statement 2.

4. Statement 4 follows from statement 1 and 2. $\Box $

$\noindent $

Next, we  find the stabilizers of the generic  ideals (see
Corollary \ref{c15Apr9}). First, we recall the definition of the
{\em wreath product} $A\wr B$ of finite groups $A$ and $B$. The
set $\Fun (B,A)$ of all functions $f: B\ra A$ is a group: $(fg)
(b) := f(b) g(b)$ for all $b\in B$ where $g\in \Fun (B,A)$. There
is a group homomorphism
$$ B\ra \Aut (\Fun (B,A)), \; b_1\mapsto (f\mapsto b_1(f):b\mapsto
f(b_1^{-1}b)).$$ Then the semidirect product $\Fun (B,A) \rtimes
B$ is called the {\em wreath product} of the groups $A$ and $B$
denoted $A\wr B$, and so the product in $A\wr B$ is given by the
rule:
$$f_1b_1\cdot f_2b_2= f_1b_1(f_2) b_1b_2, \;\; {\rm where}\;\;
f_1, f_2\in \Fun (B,A) , \;\; b_1,b_2\in B.$$ Recall that each
nonzero prime ideal $\gp$ of the algebra $\mI_n$ is a unique sum
$\gp = \sum_{i\in I} \gp_i$ of height one prime ideals. The set
$\Supp (\gp ):= \{ \gp_i\, | \, i\in I\}$ is called the {\em
support} of $\gp$.

$\noindent $

{\it Definition}. We say that a proper  ideal $\ga$ of $\mI_n$ is
{\em generic} if $\Supp (\gp ) \cap \Supp (\gq )=\emptyset$ for
all $\gp , \gq \in \Min (\ga )$ such that $\gp \neq \gq$.

\begin{corollary}\label{c15Apr9}
Let $\ga$ be a generic  ideal of the algebra $\mI_n$. The set
$\Min (\ga )$ of minimal primes over $\ga$ is the disjoint union
of its non-empty subsets, $\Min_{h_1}(\ga ) \bigcup \cdots \bigcup
\Min_{h_t}(\ga )$, where $1\leq h_1<\cdots < h_t\leq n$ and the
set $\Min_{h_i}(\ga )$ contains all the minimal primes over $\ga$
of height $h_i$. Let $n_i:= |\Min_{h_i}(\ga )|$. Then $
\St_{G_n}(\ga )= (S_m\times \prod_{i=1}^t(S_{h_i}\wr
S_{n_i}))\ltimes \mT^n \ltimes \Inn (\mI_n) $ where $m=
n-\sum_{i=1}^t n_ih_i$.
\end{corollary}

{\it Proof}. Suppose that $\Min (\ga ) = \{ \gq_1, \ldots ,
\gq_s\}$ and the sets $I_1, \ldots , I_s$ are defined in Theorem
\ref{15Apr9}. Since the ideal $\ga$ is generic, the sets $I_1,
\ldots , I_s$ are disjoint.  By Theorem \ref{15Apr9}, we have to
show that 
\begin{equation}\label{Smm}
\St_{S_m}(\{ I_1, \ldots , I_s\}) \simeq S_m\times
\prod_{i=1}^t(S_{h_i}\wr S_{n_i}).
\end{equation}
The ideal $\ga$ is generic, and so the set $\{ 1, \ldots , n\}$ is
the disjoint union $\bigcup_{i=0}^t M_i$ of its subsets where
$M_i:= \bigcup_{|I_j|=h_i}I_j$, $i=1, \ldots , t$, and $M_0$ is
the complement of the set $\bigcup_{i=1}^tM_i$. Let $S(M_i)$ be
the symmetric group corresponding to the set $M_i$ (i.e. the set
of all bijections $M_i\ra M_i$). Then each element $\s \in
\St_{G_n}(\{ I_1, \ldots , I_s\} )$ is a unique product $\s =
\s_0\s_1\cdots \s_t$ where $\s_i\in S(M_i)$. Moreover, $\s_0$ can
be an arbitrary element of $S(M_0) \simeq S_m$, and, for $i\neq
0$, the element $\s_i$ permutes the sets $\{ I_j\, | \,
|I_j|=h_i\}$ and simultaneously permutes the elements inside each
of the sets $I_j$, i.e. $\s_i\in S_{h_i}\wr S_{n_i}$. Now,
(\ref{Smm}) is obvious. $\Box $

\begin{corollary}\label{e15Apr9}
For each number $s=1, \ldots, n$, let
$\gb_s:=\prod_{|I|=s}(\sum_{i\in I} \gp_i)$ where $I$ runs through
all the subsets of the  set  $\{ 1, \ldots , n\}$ that contain
exactly $s$ elements.  The ideals $\gb_s$ are the only proper,
$\rG_n$-invariant ideals of the algebra $\mI_n$, and so there are
precisely $n+2$ $\rG_n$-invariant ideals of the algebra $\mI_n$.
\end{corollary}

{\it Proof}. By Theorem \ref{15Apr9}, the ideals  $\gb_s$ are
$\rG_n$-invariant, and they are proper. The converse follows at
once from the classification of  ideals for the algebra $\mI_n$
(Theorem \ref{a15Apr9}) and Theorem \ref{15Apr9}. The ideals
$\gb_s$ are distinct, by Theorem \ref{a15Apr9}, and so there are
precisely $n+2$ $\rG_n$-invariant ideals of the algebra $\mI_n$.
$\Box $

$\noindent $

Department of Pure Mathematics

University of Sheffield

Hicks Building

Sheffield S3 7RH

UK

email: v.bavula@sheffield.ac.uk

\end{document}